\documentclass[11pt,amsfonts]{article}
\usepackage{graphicx}
\usepackage{latexsym}
\usepackage{amssymb}
\usepackage{amsmath}
\usepackage{layout}
\newtheorem{prop}{Proposition}
\newtheorem{lemma}{Lemma}
\newtheorem{definition}{Definition}
\newtheorem{corollary}{Corollary}
\newtheorem{theorem}{Theorem}
\newtheorem{remark}{Remark}
\newtheorem{example}{Example}

\def\real{{\mathord{{\rm I\kern-2.8pt R}}}}        
\def\inte{{\mathord{{\rm I\kern-2.8pt N}}}}

\def\sZZ{{\rm Z\kern-2.8ptem{}Z}}

\def\z{{\mathchoice
  {\sZZ}
  {\sZZ}
  {\rm Z\kern-0.30em{}Z}
  {\rm Z\kern-0.25em{}Z} }}
\def\sQQ{{\kern 0.27em \vrule height1.45ex width0.03em depth0em
          \kern-0.30em \rm Q}}
\def\qu{{\mathchoice
    {\sQQ}
    {\sQQ}
  {\kern 0.225em \vrule height1.05ex width0.025em depth0em \kern-0.25em \rm Q}
  {\kern 0.180em \vrule height0.78ex width0.020em depth0em \kern-0.20em \rm Q}
        }}
\def\sCC{{\kern 0.27em \vrule height1.45ex width0.03em depth0em
          \kern-0.30em \rm C}}
\def\complex{{\mathchoice
    {\sCC}
    {\sCC}
  {\kern 0.225em \vrule height1.05ex width0.025em depth0em \kern-0.25em \rm C}
  {\kern 0.180em \vrule height0.78ex width0.020em depth0em \kern-0.20em \rm C}
        }}

\newcommand{\ba}{\begin{array}}
\newcommand{\ea}{\end{array}}
\newcommand{\be}{\begin{equation}}
\newcommand{\ee}{\end{equation}}
\newcommand{\bea}{\begin{eqnarray}}
\newcommand{\eea}{\end{eqnarray}}
\newcommand{\beaa}{\begin{eqnarray*}}
\newcommand{\eeaa}{\end{eqnarray*}}

\def\z{\zeta}

\def\n{\nu}

\font\tenmath=msbm10 \font\sevenmath=msbm7 \font\fivemath=msbm5
\newfam\mathfam \textfont\mathfam=\tenmath
\scriptfont\mathfam=\sevenmath \scriptscriptfont\mathfam=\fivemath

\def \={{\buildrel {\rm (law)} \over =}}

\def\qed{ \hfill \vrule width.25cm height.25cm depth0cm\smallskip}

\newcommand{\basa}{\begin{assumption}}
\newcommand{\easa}{\end{assumption}}

\newcommand{\bas}{\begin{assum}}
\newcommand{\eas}{\end{assum}}

\newcommand{\ignore}[1]{}
\begin{document}

\renewcommand{\thefootnote}{\fnsymbol{footnote}}

\title{Analysis of the Rosenblatt process}

\vskip1cm

\author{Ciprian A.
Tudor \\
SAMOS/MATISSE, Centre d'Economie de
La Sorbonne,\\ Universit\'e de Panth\'eon-Sorbonne Paris 1,\\
90,
rue de Tolbiac, 75634 Paris Cedex 13, France.}
\maketitle

 \begin{abstract}
We analyze {\em the Rosenblatt process } which is a selfsimilar
process with stationary increments and which   appears as limit in
the so-called {\em Non Central Limit Theorem } (Dobrushin and Major
(1979), Taqqu (1979)). This process is non-Gaussian and it lives in
the second Wiener chaos. We give its representation as a
Wiener-It\^o multiple integral with respect to the Brownian motion
on a finite interval
 and we develop a stochastic calculus with respect to it by using both pathwise type calculus and Malliavin calculus.
 \end{abstract}

\vskip0.3cm

{\bf  2000 AMS Classification Numbers: }60G12, 60G15, 60H05, 60H07

 \vskip0.3cm

{\bf Key words: } Non Central Limit Theorem, Rosenblatt process,
fractional Brownian motion, stochastic calculus via regularization,
Malliavin calculus, Skorohod integral.

\vskip0.3cm

 \section{Introduction}
A selfsimilar object is exactly of approximately similar to a part
of itself. Selfsimilar processes are invariant in distribution under
suitable scaling. They are of considerable interest in practice
since aspects of the selfsimilarity appear in different phenomena
like telecommunications, economics, hydrology or turbulence. We
refer to the work of Taqqu \cite{Taqqu86} for a guide on the
appearance  of the selfsimilarity in many applications and to the
monographs by Samorodnitsky and Taqqu \cite{ST} and by Embrechts and
Maejima \cite{EM} for complete expositions  on selfsimilar
processes.

In this work we analyze a special class of selfsimilar processes
that are limits in the so called {\em Non Central Limit Theorem  }
(see Dobrushin and Major \cite{DM} or Taqqu \cite{Taqqu79}).  Let us
briefly recall the general context.

Consider $(\xi _{n}) _{n\in \mathbb{Z}}$ a stationary Gaussian
sequence with mean zero and variance $1$ such that its correlation
function satisfies
\begin{equation}
\label{corr} r(n) := \mathbf{E}\left( \xi _{0} \xi _{n} \right) =
n^{\frac{2H-2} {k}} L(n)
\end{equation}
with $H\in (\frac{1}{2}, 1) $  and $L$ is a slowly varying function
at infinity (see e.g. \cite{EM}). Denote by $H_{m}(x)$ the Hermite polynomial of degree
$m$ given by $ H_{m}(x)=
(-1)^{m}e^{\frac{x^{2}}{2}}\frac{d^{m}}{dx^{m}}e^{-\frac{x^{2}}{2}}
$ .  Let $g $ be a function such that $\mathbf{E}(g(\xi _{0} ))=0$ and
$\mathbf{E}(g(\xi _{0}^{2}))<\infty. $ Suppose that $g$ {\em has Hermite rank
} equal to $k$; that is, if $g$ admits the following expansion in
Hermite polynomials
\begin{equation*}
g(x)= \sum _{j\geq 0} c_{j} H_{j}(x), \hskip0.5cm c_{j}=
\frac{1}{j!} \mathbf{E}\left( g(\xi _{0} H_{j} (\xi _{0})) \right)
\end{equation*}
then
\begin{equation*}
k=\min \{ j; c_{j}\not= 0\}.
\end{equation*}
 Since $\mathbf{E}\left[ g(\xi
 _{0})\right] =0$,  we have $k\geq 1$. Then the Non Central Limit
 Theorem (\cite{DM}, \cite{Taqqu79}) says that
 \begin{equation*}
 \frac{1}{n^{H}} \sum _{ j=1} ^{[nt]} g(\xi_{j})
 \end{equation*}
 converges as $n \to \infty $ in the sense of finite dimensional
 distributions  to the process
 \begin{equation}
\label{hermite} Z^{k}_{H} (t)= c(H,k)
\int_{\mathbb{R}^{k}}\int_{0}^{t} \left( \prod _{j=1}^{k} (s-y_{i}
)_{+} ^{ -\left( \frac{1}{2} + \frac{1-H}{k} \right) }  \right) ds
dB(y_{1} )\ldots dB(y_{k}),
\end{equation}
where $x_{+}=\max(x,0)$ and the above integral is a multiple
Wiener-It\^o stochastic integral with respect to a Brownian motion
$B(y))_{y\in \mathbb{R}}$ (see \cite{N} for the definition). The
constant $c(H,k)$ is positive and it will be taken such that
$\mathbf{E}\left(  Z^{k}_{H} (1) ^{2} \right)=1$. The process
$(Z_{H}^{k}(t) )_{t\geq 0} $ is called {\em the Hermite process }
and it is $H$-selfsimilar in the sense that for any $c>0$,
$(Z^{k}_{H}(ct)) =  ^{(d)} (c^{H} Z^{k}_{H} (t))$, where $ " =^{(d)}
"$ means equivalence of all finite dimensional distributions, and it
has stationary increments.

When $k=1$ the process given by (\ref{hermite}) is nothing else that
the {\em fractional Brownian motion (fBm) } with Hurst parameter
$H\in ( \frac{1}{2}, 1)$. For $k\geq 2$ the process is not Gaussian.
If $k=2$ then the process (\ref{hermite}) is known as {\em the
Rosenblatt process} (it has actually called in this  way by M. Taqqu
in \cite{Taqqu75}).

The fractional Brownian motion is of course the most studied process in the
class of Hermite processes due to its significant importance in
modeling. A stochastic calculus with respect to it has  been
intensively developed in the last decade. We refer, among others, to
\cite{AMN}, \cite{AN}, \cite{DU}, \cite{GNRV}.

Our main interest consists here in the study,  from the stochastic
calculus point of view, of the Rosenblatt process. Although it
received a less important attention than the fractional Brownian
motion, this process is still of interest in practical  applications because
of its self-similarity, stationarity of increments and long-range
dependence. There exists a consistent literature that focuses on
different theoretical aspects of the Rosenblatt processes. Let us
recall some of these works. For example, extremal properties of the
Rosenblatt distribution have been studied by J.M.  Albin in \cite{A1}
and \cite{A2}. The rate of convergence to the Rosenblatt process in
the Non Central Limit Theorem has been given by Leonenko and Ahn
\cite{LA}. Pipiras \cite{Pi} and Pipiras and Abry  \cite{PA}  studied the wavelet-type expansion of
the Rosenblatt process. A law of iterated logarithm  has been given
is \cite{Kuelb}.

Among the applications of the Rosenblatt process in statistics or
econometrics, we mention the following.
\begin{description}
\item{$\bullet$ } In the unit root testing problem  with errors being nonlinear
transforms of linear processes with long-range dependence, the
asymptotic distributions in the model  are shown in \cite{Wu} to be
functionals of Hermite processes.
\item{$\bullet$ } limiting distributions of the parabolically
rescaled solutions of the heat equation with singular non-Gaussian
data have similar behavior to the Rosenblatt distribution (see
\cite{LW})
\item{$\bullet$ } the Rosenblatt distribution also appears to be the
asymptotic distribution of an estimator related to the
semiparametric bootstrap approach to hypothesis tests (see
\cite{Hall}) or to the estimation of the long-range dependence parameter (\cite{KG})
\end{description}
Besides these more or less practical applications of the Rosenblatt
process, denoted in the following by $Z$, our motivation is also
theoretical; it comes from the recent intensive interest to push
further the stochastic calculus with respect to more and more
general integrator processes. We believe that this process
constitutes an interesting and instructive example where the recent
developed techniques of the generalized stochastic calculus can find
a significant test bench.

We actually  use the two principal methods to develop
 a stochastic integration theory: the pathwise type calculus and the Malliavin calculus/Skorohod integration.
 The first approach (that includes essentially the rough paths analysis,
  see \cite{QL}, and the stochastic calculus via regularization, see \cite{RV1})
  can be directly applied to the Rosenblatt process because of its regular paths and of the nice covariance structure;
  a pathwise It\^o formula can be written and Stratonovich stochastic
  equation with $Z$ as noise can be considered. The Malliavin calculus and the Skorohod integration are in general
  connected in a deeper way to the Gaussian structure of the integrator process and as it will be seen in the present work, Skorohod It\^o formula can be derived only in particular cases. Although the formula we obtain is rather complicated and not easily tractable, the principal signification of the result is the fact that one can precisely see here that the Gaussian nature of the integrator process is decisive in the stochastic integration theory; once we go out from the Gaussian context, one cannot obtain It\^o's formulas that end by a second derivative term.

\smallskip

We organized our paper as follows. Section 2 presents  basic properties  of the Rosenblatt process. In particular
we prove a stochastic integral representation on a finite integral that will be useful
 for the construction of the stochastic calculus.
 In Section 3 we introduce Wiener integrals with respect to $Z$ by following the ideas in \cite{MaTu} and \cite{KRT}.
 We define in Section 4 the Hilbert-valued Rosenblatt process and we consider stochastic evolution equations with this process as noise.
 Section 5 describes the application of the
  stochastic calculus via regularization introduced by F. Russo and P. Vallois in \cite{RV1} to the Rosenblatt process and
   in Section 6 we  discuss the Skorohod (divergence) integral: we define the integral and
   we give conditions that ensure the integrability and the continuity of the indefinite integral process.
In Section 7 we prove the relation between the pathwise and the divergence integrals:
 here the pathwise integral is equal to the Skorohod integral plus two trace terms
  (in the fBm case there is  only a trace term). Finally Section 8 contains a discussion on the It\^o formula  in the Skorohod sense.

  \vskip0.5cm

 \section{On the Rosenblatt process}
 In this section we will analysis  some basic properties of the
 Rosenblatt process; in particular we are interested in its representation as a stochastic integral
 on a finite interval.  As we said, this  process is obtained
 by taking $k=2$ in the relation (\ref{hermite}), so
 \begin{equation}
 \label{rose1}
Z_{2}(t):=Z(t) = a(H) \int_{\mathbb{R}} \int _{\mathbb{R}} \left(
\int_{0}^{t} (s-y_{1}) _{+}^{-\frac{2-H}{2}}(s-y_{2})_{+} ^{-\frac{2-H}{2}}ds \right)
 dB(y_{1})dB(y_{2})
 \end{equation}
 where $ (B(y), y\in \mathbb{R})$  is a standard Brownian motion on
 $\mathbb{R}$. The constant $a(H) $ is a positive normalizing constant and it is chosen such that
 $\mathbf{E}(Z(1)^{2})=1$. It follows actually from \cite{MaTu} that
 \begin{equation*}
a(H)^{2}= \left( \frac{\beta (\frac{H}{2}, H-1) ^{2}}{2H(2H-1)}
\right) ^{-1}
 \end{equation*}Recall that the process $(Z(t))_{ t\in [0,T]}$ is
 selfsimilar of order $H$ and it has stationary increments; it admits a
 H\"older continuous version of order $\delta <H$. Since $H\in
 (\frac{1}{2}, 1)$, it follows that the process $Z$ exhibits
 long-range dependence.

 \vskip0.3cm

 Since our main interest consists in the construction of the
 stochastic calculus with respect to the process $Z$, the
 representation (\ref{rose1}) is not very convenient; as in the fBm
 case, we would  like to represent $Z_{t}$ as a stochastic integral with
 respect to a Brownian motion with time interval $[0,T]$. Recall
 that the fBm with $H>\frac{1}{2}$ can be written as
 \begin{equation}
 \label{B1}
B_{t}^{H}= \int_{0}^{t} K^{H}(t,s)dW_{s}
 \end{equation}
 with $(W_{t}, t\in [0,T])$ a standard Wiener process and
 \begin{equation}
 \label{K}
K^{H}(t,s)= c_{H} s^{\frac{1}{2}-H} \int _{s}^{t}
(u-s)^{H-\frac{3}{2}} u^{H-\frac{1}{2}}  du
 \end{equation}
 where $t>s$ and
 \begin{equation}
 \label{ch}
c_{H} =\left( \frac{ H(2H-1) }{\beta( 2-2H, H-\frac{1}{2}) } \right)
^{\frac{1}{2}}.
\end{equation}
Note that to prove the representation (\ref{B1}) (at least in law)
it suffices to see that the right member has  the
same covariance $R$ as the fBm; otherwise, it can be easily seen from  the expression
of the kernel $K$ that the right member in (\ref{B1}) is
$H$-selfsimilar with stationary increments and as a consequence it
cannot be nothing else but a fractional Brownian motion with
parameter $H$.

Since  the Rosenblatt process is not Gaussian,   the proof in its
case of a similar representation to (\ref{B1})  needs a
supplementary argument; in fact we have the following

\begin{prop}
\label{egald} Let $K$ be the kernels (\ref{K}) and let $(Z(t))
_{t\in [0,T]}$ be a Rosenblatt process with parameter $H$. Then it
holds that
\begin{equation}
\label{rose2} Z(t)=^{(d)} d(H) \int _{0}^{t}\int_{0}^{t} \left[
\int_{ y_{1} \vee y_{2} }^{t} \frac{\partial K^{H'} }{\partial u}
(u,y_{1} )  \frac{\partial K^{H'} }{\partial u} (u,y_{2} )du
\right]dB(y_{1}) dB(y_{2})
\end{equation}
where $(B_{t}, t \in [0,T])$ is a Brownian motion,
\begin{equation}
\label{H'} H'=\frac{H+1}{2}
\end{equation}
and \begin{equation} \label{dh} d(H)= \frac{1}{H+1} \left(
\frac{H}{2(2H-1)} \right) ^{-\frac{1}{2}}.
\end{equation}
\end{prop}

\begin{remark}
i) The constant $d(H)$ is a normalizing constant, it has been chosen
such that $\mathbf{E}(Z(t)Z(s)) = \frac{1}{2}\left( t^{2H} + s^{2H } -\vert
t-s \vert ^{2H} \right)$. Indeed,
\begin{eqnarray*}
\mathbf{E}(Z(t)Z(s))&=& 2d(H) ^{2} \int _{0}^{t\wedge s}\int_{0}^{t\wedge s} dy_{1}dy_{2}\\
&&\times  \left(   \int_{ y_{1} \vee y_{2} }^{t} \int_{ y_{1} \vee
y_{2} }^{s}\frac{\partial K^{H'} }{\partial u} (u,y_{1} )
\frac{\partial K^{H'} }{\partial u} (u,y_{2} )\frac{\partial K^{H'}
}{\partial u} (v,y_{1} ) \frac{\partial K^{H'} }{\partial v}
(v,y_{2} ) dudv
\right)\\
&=&2d(H) ^{2}\int _{0}^{t}\int_{0}^{s}dvdu \left( \int_{0}^{u\wedge
v } \frac{\partial K^{H'} }{\partial u} (u,y_{1} )\frac{\partial
K^{H'} }{\partial u} (v,y_{1} )dy_{1}\right) ^{2}\\
&=&2d(H) ^{2} (H'(2H'-1) )^{2}\int _{0}^{t}\int_{0}^{s} \vert
u-v\vert ^{2H-2}dvdu  = R(t,s).
\end{eqnarray*}

\vskip0.5cm

 ii) It can be seen without without difficulty that the process
\begin{equation*}
Z'(t):=d(H) \int _{0}^{t}\int_{0}^{t} \left[ \int_{ y_{1} \vee y_{2}
}^{t} \frac{\partial K^{H'} }{\partial u} (u,y_{1} ) \frac{\partial
K^{H'} }{\partial u} (u,y_{2} )du \right]dB(y_{1}) dB(y_{2})
\end{equation*}defines a $H$ selsimilar process  with stationary
increments. Indeed, for any $c>0$,
\begin{eqnarray*}
Z'(ct)&=& \int_{0}^{ct} \int _{0}^{ct} \left[\int_{ y_{1} \vee y_{2}
}^{ct} \frac{\partial K^{H'} }{\partial u} (u,y_{1} ) \frac{\partial
K^{H'} }{\partial u} (u,y_{2} )du \right]dB(y_{1})
dB(y_{2})\\
 &=&\int_{0}^{ct}
\int _{0}^{ct}\left[ \int _{\frac{y_{1}}{c}\vee
\frac{y_{2}}{c}}^{t}\frac{\partial K^{H'} }{\partial u} (cu,y_{1} )
\frac{\partial K^{H'} }{\partial u} (cu,y_{2} )c du\right]dB(y_{1})
dB(y_{2})\\
&=&\int _{0}^{t}\int_{0}^{t} \left[ \int_{ y_{1} \vee y_{2} }^{t}
\frac{\partial K^{H'} }{\partial u} (cu,cy_{1} ) \frac{\partial
K^{H'} }{\partial u} (cu,cy_{2} )cdu \right]dB(cy_{1}) dB(cy_{2})
\end{eqnarray*}
and since $B(cy)= ^{(d)} c^{\frac{1}{2} }B(y) $ and $\frac{\partial
K^{H'} }{\partial u} (cu,cy_{i} ) =
c^{H^{'}-\frac{3}{2}}\frac{\partial K^{H'} }{\partial u} (u,y_{i} )$
we obtain $Z(ct)=^{(d)} c^{H} Z(t)$.

The fact that $Z'$ has stationary increments follows from the
relation \begin{equation*} K^{H'}(t+h,s)-K^{H'}(t,s)=
K^{H'}(t-s,h)\end{equation*} for any $s,t\in [0,T]$, $s<t$ and
$h>0$.
\end{remark}

\vskip0.2cm

{\bf Proof of Proposition \ref{egald}: } Let us denote by $Z'(t)$
the right hand side of (\ref{rose2}). Consider $b_{1}, \ldots ,
b_{n} \in \mathbb{R}$ and $t_{1}, \ldots , t_{n}\in [0,T]$. We need
to show that the random variables
\begin{equation*}
\sum_{l=1}^{n} b_{l} Z(t_{l}), \hskip0.5cm \sum_{l=1}^{n} b_{l}
Z'(t_{l})
\end{equation*}
have the same distribution.

\vskip0.3cm

We will use the following criterium by Fox and Taqqu (see
\cite{FT}): If $f\in L^{2}([0,T]^{2})$ is a symmetric function, then
the law of the multiple Wiener-It\^o integral $I_{2}(f)$ is uniquely
determined by its {\it cumulants}, where the $m$th cumulant of $f$
is given by
\begin{equation}
\label{cum} c_{m}(f)=\frac{(m-1)! } {2} 2^{m} \int _{\mathbb{R}
^{m}}f(x_{1}, x_{2} ) f(x_{2}, x_{3}) \ldots
f(x_{m-1},x_{m})f(x_{m},x_{1}) dx_{1}\ldots dx_{m} .
\end{equation}
In other words, if two symmetric functions $f,g\in L^{2}([0,T]^{2})$
have the same cumulants, then the multiple Wiener-It\^o  integrals
of order two  $I_{2}(f) $ and $I_{2}(g)$ have the same law.

\vskip0.5cm

We will show that, for every $t,s\in [0,T]$, the random variables
$Z_{t}+Z_{s} $ and $Z'_{t}+ Z'_{s} $ have the same law; the general
case case will follow by a similar calculation. It holds that
\begin{equation*}
Z'_{t}+ Z'_{s} = I_{2}\left( f_{t,s} \right)
\end{equation*}
where
\begin{eqnarray}
 f_{t,s}(y_{1}, y_{2})&=&
1_{[0,t]}(y_{1})1_{[0,t]}(y_{2})\int _{y_{1}\vee y_{2}}
^{t}\frac{\partial K^{H'} }{\partial u} (u,y_{1} ) \frac{\partial
K^{H'} }{\partial u} (u,y_{2} )du\nonumber \\
&&+1_{[0,s]}(y_{1})1_{[0,s]}(y_{2})\int _{y_{1}\vee y_{2}}
^{s}\frac{\partial K^{H'} }{\partial u} (u,y_{1} ) \frac{\partial
K^{H'} }{\partial u} (u,y_{2} )du v
\end{eqnarray}
We have denoting by $a_{m}:=\frac{(m-1)! } {2} 2^{m}d(H)^{m}  $,
\begin{eqnarray*}
&&c _{m}(f_{s,t})\\&=& a(m) \int _{\mathbb{R}^{m}}f_{t,s}(y_{1}, y_{2})
\ldots f_{t,s}(y_{m}, y_{1}) dy_{1}\ldots dy_{m} \\
&=& a(m) \int _{\mathbb{R}^{m}}dy_{1}\ldots dy_{m} \\
&&\times \left( \int _{y_{1}\vee y_{2}} ^{t}\frac{\partial K^{H'}
}{\partial u} (u_{1},y_{1} ) \frac{\partial K^{H'} }{\partial u}
(u_{1},y_{2} )du_{1}+\int _{y_{1}\vee y_{2}} ^{s}\frac{\partial
K^{H'} }{\partial u} (u_{1},y_{1} ) \frac{\partial K^{H'} }{\partial
u} (u_{1},y_{2} )du_{1}\right) \\
&&\times \left( \int _{y_{2}\vee y_{3}} ^{t}\frac{\partial K^{H'}
}{\partial u} (u_{2},y_{2} ) \frac{\partial K^{H'} }{\partial u}
(u_{2},y_{3} )du_{2}+\int _{y_{2}\vee y_{3}} ^{s}\frac{\partial
K^{H'} }{\partial u} (u_{2},y_{2} ) \frac{\partial K^{H'} }{\partial
u} (u_{2},y_{3} )du_{2}\right)\\
&&\times \ldots   \\
&&\times \left( \int _{y_{m}\vee y_{1}} ^{t}\frac{\partial K^{H'}
}{\partial u} (u_{m},y_{m} ) \frac{\partial K^{H'} }{\partial u}
(u_{m},y_{1} )du_{m}+\int _{y_{m}\vee y_{1}} ^{s}\frac{\partial
K^{H'} }{\partial u} (u_{m},y_{1} ) \frac{\partial K^{H'} }{\partial
u} (u_{m},y_{m} )du_{m}\right)
\end{eqnarray*}
and by classical Fubini theorem
\begin{eqnarray}
&&c_{m}(f_{s,t})\nonumber\\
&=& a(m)\sum_{t_{j} \in \{t,s\}} \int_{0}^{t_{1}} \ldots \int_{0}^{t_{m}} du_{1}\ldots du_{m} \nonumber \\
&&\times \left( \int _{0} ^{ u_{1} \wedge u_{m}}\frac{\partial
K^{H'} }{\partial u_{1}} (u_{1},y_{1} )\frac{\partial K^{H'}
}{\partial u_{m}} (u_{m},y_{1} )dy_{1}\right) \nonumber  \\
&&\times \left( \int _{0} ^{ u_{1} \wedge u_{2}} \frac{\partial
K^{H'} }{\partial u_{1}} (u_{1},y_{2} )\frac{\partial K^{H'}
}{\partial u_{2}} (u_{2},y_{2} )dy_{2} \right) \nonumber\\
&&\ldots \nonumber  \\
 &&\times \int_{0} ^{u_{m-1} \wedge u_{m}}\frac{\partial K^{H'}
}{\partial u_{m-1}} (u_{m},y_{m} )\frac{\partial K^{H'} }{\partial
u_{m}} (u_{m},y_{m} )dy_{m}\nonumber \\
&=& a(m)  \sum_{t_{j} \in \{t,s\}} \int_{0}^{t_{1}} \ldots
\int_{0}^{t_{m}} du_{1} \ldots du_{m}\nonumber \\
&&   \left|  u_{1}-u_{2}\right| ^{2H'-2} \left|
u_{2}-u_{3}\right| ^{2H'-2}\ldots  \left| u_{m}-u_{1}\right|
^{2H'-2}.\label{cf}
\end{eqnarray}
with $a(m)=a(m)\left( H'(2H'-1)\right) ^{m}$.

 The computation of the cumulant of $Z_{t}+Z_{s}$ is similar. Indeed, we can write, for $s,t\in [0,T]$,
\begin{equation*}
Z(t)+Z(s)= I_{2} (g_{s,t} )
\end{equation*}
where
\begin{equation*}
g_{s,t}= a(H) \left( \int _{0}^{t}  (u-y_{1})_{+} ^{\frac{H-2} {2}}
(u-y_{2 }) _{+} ^{\frac{H-2} {2}} du    +\int _{0}^{s}  (u-y_{1})_{+}
^{\frac{H-2} {2}} (u-y_{2 })_{+}  ^{\frac{H-2} {2}} du \right)
\end{equation*}
and the $m$th cumulant of the kernel $g_{s,t}$ is given by
\begin{eqnarray*}
c_{m}(g_{s,t}) &= & b(m) \int _{\mathbb{R} ^{m}}dy_{1}\ldots
dy_{m}\\
&&\left( \int _{0}^{t}  (u_{1}-y_{1})_{+} ^{\frac{H-2} {2}}
(u_{1}-y_{2 }) _{+} ^{\frac{H-2} {2}} du_{1}    +\int _{0}^{s}
(u_{1}-y_{1}) _{+} ^{\frac{H-2}
{2}} (u_{1}-y_{2 })_{+}  ^{\frac{H-2} {2}} du_{1} \right) \\
&&\left( \int _{0}^{t}  (u_{2}-y_{2})_{+} ^{\frac{H-2} {2}}
(u_{2}-y_{3 }) _{+} ^{\frac{H-2} {2}} du_{2}    +\int _{0}^{s}
(u_{2}-y_{2})_{+}  ^{\frac{H-2} {2}} (u_{2}-y_{3 })_{+}  ^{\frac{H-2} {2}}
du_{2} \right)\\
&&\ldots  \\
&&\left( \int _{0}^{t}  (u_{m}-y_{m})_{+} ^{\frac{H-2} {2}}
(u_{m}-y_{1 }) _{+} ^{\frac{H-2} {2}} du_{1}    +\int _{0}^{s}
(u_{m}-y_{m}) _{+} ^{\frac{H-2} {2}} (u_{m}-y_{1 }) _{+} ^{\frac{H-2} {2}}
du_{m} \right)\\
&=& b(m)\sum_{t_{j} \in \{t,s\}} \int_{0}^{t_{1}} \ldots \int_{0}^{t_{m}} du_{1}\ldots du_{m} \\
&=& \int _{\mathbb{R}}(u_{1}-y_{1})_{+} ^{\frac{H-2}
{2}}(u_{m}-y_{1})_{+} ^{\frac{H-2} {2}}dy_{1} \int _{\mathbb{R}}
(u_{1}-y_{2})_{+} ^{\frac{H-2} {2}}(u_{2}-y_{2})_{+} ^{\frac{H-2}
{2}}dy_{2}\\
&& \ldots  \ldots \int_{\mathbb{R}}(u_{m-1}-y_{m})_{+} ^{\frac{H-2}
{2}}(u_{m}-y_{m})_{+} ^{\frac{H-2} {2}}dy_{m}.
\end{eqnarray*}
Since for any $a>0$
\begin{equation*}
 \int_{\mathbb{R}} (u-y)_{+} ^{a-1} (v-y)_{+}^{a-1}
dy =\beta (a, 2a-1) \vert u-v\vert ^{2a-1}
\end{equation*}
we get
\begin{eqnarray}
 c_{m}(g_{s,t})& =& b(m) \beta (\frac{H}{2}, H-1)
^{m}\sum_{t_{j} \in \{t,s\}} \int_{0}^{t_{1}} \ldots
\int_{0}^{t_{m}}  du_{1} \ldots
du_{m}\nonumber \\
&&\left| u_{1}-u_{2}\right| ^{2H'-2} \left| u_{2}-u_{3}\right|
^{2H'-2}\ldots \left| u_{m}-u_{1}\right| ^{2H'-2}\label{cg}
\end{eqnarray}
and it remains to observe that $a'(m)= b(m) $ which implies that
(\ref{cf}) equals (\ref{cg}). \qed

\vskip0.5cm

{\em From now on we will use the version of the Rosenblatt  process
given by the right side of  (\ref{rose2}). }

\vskip0.5cm

We will finish this section  by proving that the Rosenblatt process
possesses a similar property to the fBm, that is, it can be
approximated by a sequence of semimartingales (here actually, since
$H>\frac{1}{2}$, by a sequence of bounded variation processes).
 In the fBm case, the property is
inherited by the divergence integral (see \cite{AMN}, \cite{Bia}, \cite{Mi}); this fact
can be used  to construct financial models with the Rosenblatt
process as noise (see \cite{Mi}).

The basic observation is that, if one interchanges formally the
stochastic and Lebesque integrals in (\ref{rose2}), one gets
\begin{equation*}
Z(t)"=" \int_{0}^{t} \left( \int_{0}^{u} \int _{0}^{u}\frac{\partial
K^{H'} }{\partial u} (u,y_{1})\frac{\partial K^{H'} }{\partial u}
(u,y_{2})dB(y_{1})dB(y_{2}) \right) du
\end{equation*}
but the above expression cannot hold because the kernel
$\frac{\partial K^{H'} }{\partial u} (u,y_{1})\frac{\partial K^{H'}
}{\partial u} (u,y_{2}) $ does not belong to $L^{2}([0,T]^{2})$ since the partial derivative $\frac{\partial K^{H'} }{\partial u} (u,y_{1})$ behaves on the diagonal as $(u-y_{1}) ^{\frac{H-2}{2}}$.

Let us define, for every $\varepsilon >0$,
\begin{eqnarray*}
Z^{\varepsilon} (t)&=& d(H) \int _{0}^{t}\int_{0}^{t} \left[ \int_{
y_{1} \vee y_{2} }^{t} \frac{\partial K^{H'} }{\partial u}
(u+\varepsilon ,y_{1} )  \frac{\partial K^{H'} }{\partial u}
(u+\varepsilon,y_{2} )du \right]dB(y_{1}) dB(y_{2})\\
&=& \int_{0}^{t} \left(  \int_{0}^{u} \int _{0}^{u}\frac{\partial
K^{H'} }{\partial u} (u+\varepsilon,y_{1})\frac{\partial K^{H'}
}{\partial u} (u+\varepsilon,y_{2})dB(y_{1})dB(y_{2}) \right) du\\
&:=& \int_{0}^{t} A_{\varepsilon }(u)du.
\end{eqnarray*}
Since $A_{\varepsilon} \in L^{2}([0,T]\times \Omega )$ for every
$\varepsilon >0$ and it is adapted, it follows that the process $Z^{\varepsilon}$ is a
semimartingale.

\begin{prop}
For every $t\in [0,T]$, $Z^{\varepsilon}(t)\to Z(t)$ in
$L^{2}(\Omega )$.
\end{prop}
{\bf Proof: } We have
\begin{eqnarray*}
&&Z^{\varepsilon}(t)-Z(t)=\int_{0}^{t}\int_{0}^{t} dB(y_{1})dB(y_{2})\\
&&   \left( \int_{ y_{1} \vee y_{2} }^{t} \left( \frac{\partial
K^{H'} }{\partial u} (u+\varepsilon ,y_{1} ) \frac{\partial K^{H'}
}{\partial u} (u+\varepsilon,y_{2} )-\frac{\partial K^{H'}
}{\partial u} (u ,y_{1} )  \frac{\partial K^{H'} }{\partial u}
(u,y_{2} )\right) du \right)
\end{eqnarray*}
and
\begin{eqnarray*}
\mathbf{E}\left|  Z^{\varepsilon}(t)-Z(t)  \right| ^{2}
&=&2\int_{0}^{t}\int_{0}^{t} dy_{1}dy_{2} \int_{ y_{1} \vee y_{2}
}^{t} \int_{
y_{1} \vee y_{2} }^{t} dvdu\\
&&\left( \frac{\partial K^{H'} }{\partial u} (u+\varepsilon ,y_{1} )
\frac{\partial K^{H'} }{\partial u} (u+\varepsilon,y_{2}
)-\frac{\partial K^{H'} }{\partial u} (u ,y_{1} )  \frac{\partial
K^{H'} }{\partial u} (u,y_{2} )\right) \\
&& \left( \frac{\partial K^{H'} }{\partial v} (v+\varepsilon ,y_{1}
) \frac{\partial K^{H'} }{\partial v} (v+\varepsilon,y_{2}
)-\frac{\partial K^{H'} }{\partial v} (v ,y_{1} )  \frac{\partial
K^{H'} }{\partial v} (v,y_{2} )\right)
\end{eqnarray*}
Clearly the quantity $\left( \frac{\partial K^{H'} }{\partial u}
(u+\varepsilon ,y_{1} ) \frac{\partial K^{H'} }{\partial u}
(u+\varepsilon,y_{2} )-\frac{\partial K^{H'} }{\partial u} (u ,y_{1}
)  \frac{\partial K^{H'} }{\partial u} (u,y_{2} )\right) $ converges
to zero as $\varepsilon \to 0$ for every $u,y_{1}, y_{2}$  and the
conclusion follows by the dominated convergence theorem.\qed

\section{Wiener integrals}

The covariance structure of the Rosenblatt process allows to
construct Wiener integrals with respect to it. We refer to Maejima
and Tudor  \cite{MaTu} for the definition of Wiener integrals with
respect to general Hermite processes and to Kruk and al. \cite{KRT}
for a more general context. Let us recall the main points and
translate   this construction in our context.

One note that
\begin{equation*}
Z(t)= \int_{0}^{T} \int_{0}^{T} I\left( 1_{[0,t]} \right) (y_{1},
y_{2}) dB(y_{1})dB(y_{2})
\end{equation*}
where the operator $I$ is defined on the set of functions $f:[0,T]
\to \mathbb{R}$ and takes values in the set of functions
$g:[0,T]^{2}\to \mathbb{R}^{2} $ and it is given by
\begin{equation}
\label{I} I(f) (y_{1}, y_{2}) =d(H)\int _{y_{1}\vee y_{2}}^{T}
f(u)\frac{\partial K^{H'} }{\partial u} (u ,y_{1} )  \frac{\partial
K^{H'} }{\partial u} (u,y_{2} ) du.
\end{equation}
If $f$ is an element of the set ${\cal{E}}$ of step functions on
$[0,T]$ of the form
\begin{equation}
\label{step} f=\sum _{i=0}^{n-1}a_{i}1  _{(t_{i}, t_{i+1} ]},
\hskip0.5cm t_{i}\in [0,T]
\end{equation}
then we naturally define its  Wiener integral with respect to $Z$ as
\begin{equation}
\label{wiener1} \int_{0}^{T} f(u)dZ(u):= \sum _{i=0}^{n-1} a_{i}
\left( Z_{t_{i+1}}-Z_{t_{i}}\right) = \int_{0}^{T}\int_{0}^{T}
I(f)(y_{1}, y_{2}) dB(y_{1})dB(y_{2}).
\end{equation}
Let ${\cal{H}}$ be the set of functions $f$ such that
\begin{equation}
\label{normaH1} \Vert f\Vert _{{\cal{H}}}^{2}:= 2\int_{0}^{T} \int
_{0}^{T}I(f) (y_{1}, y_{2})^{2} dy_{1}dy_{2} <\infty.
\end{equation}
It can be seen that
\begin{eqnarray*}
\Vert f\Vert _{{\cal{H}}}^{2} &=& 2d(H)^{2}\int_{0}^{T} \int
_{0}^{T} \left( \int _{y_{1}\vee y_{2}}^{T} f(u)\frac{\partial
K^{H'} }{\partial u} (u ,y_{1} )  \frac{\partial K^{H'} }{\partial
u} (u,y_{2} ) du\right) ^{2} dy_{1}dy_{2} \\
&=& 2d(H)^{2}\int_{0}^{T} \int _{0}^{T}dy_{1}dy_{2}\int _{y_{1}\vee
y_{2}}^{T}
\int _{y_{1}\vee y_{2}}^{T}dvdu \\
&&  f(u)f(v)\frac{\partial K^{H'} }{\partial u} (u ,y_{1} )
\frac{\partial K^{H'} }{\partial u} (u,y_{2} ) \frac{\partial K^{H'}
}{\partial v} (v ,y_{1} ) \frac{\partial
K^{H'} }{\partial v} (v,y_{2} ) \\
&=&2d(H)^{2}\int_{0}^{T} \int _{0}^{T} \left( \int_{0}^{u\wedge v}
\frac{\partial K^{H'} }{\partial u} (u ,y_{1} ) \frac{\partial
K^{H'} }{\partial v} (v ,y_{1} )dy_{1}\right) ^{2} dvdu \\
&=& H(2H-1)\int_{0}^{T} \int _{0}^{T} f(u)f(v) \vert u-v\vert
^{2H-2} dvdu.
\end{eqnarray*}
It can be proved as in \cite{MaTu} or \cite{KRT} that the mapping
\begin{equation*}
f\to \int_{0}^{T} f(u)dZ(u)
\end{equation*}
defines an isometry from ${\cal{E}}$ to $L^{2}(\Omega)$ and it can
be extended by continuity to an isometry from ${\cal{H}}$ to
$L^{2}(\Omega)$ because ${\cal{E}}$ is dense in ${\cal{H}}$ (see
\cite{PiTa1}). We will call this extension the Wiener integral of
$f\in {\cal{H}}$  with respect to $Z$.

\begin{remark}
It follows from Pipiras and Taqqu (see \cite{PiTa1}) that the space
${\cal{H}}$ contains not only functions but its elements could be
also distributions. Therefore it is suitable to know subspaces of
${\cal{H}}$ that are spaces of functions. A such subspace is $\left|
{\cal{H}}\right| $ where
\begin{equation*}
 \left| {\cal{H}}\right| =\{
f:[0,T]\to \mathbb{R} |  \int _{0} ^{T}\int_{0} ^{T}\vert f(u)\vert
\vert f(v)\vert \vert u-v\vert ^{2H-2} dvdu <\infty \}.
\end{equation*}
It actually holds
\begin{equation*}
L^{\frac{1}{H} }([0,T]) \subset  \left| {\cal{H}}\right| \subset
{\cal{H}}.
\end{equation*}
The space $\left| {\cal{H}}\right|$ (and hence ${\cal{H}}$) is not
complete with respect to the norm $\Vert \cdot \Vert _{{\cal{H}}}$
but it is a Banach space with respect to the norm
\begin{equation*}
 \Vert f\Vert ^{2}_{\left| {\cal{H}}\right|
}=H(2H-1)\int _{0}^{T}\int_{0}^{T} \vert f(u)\vert \vert
f(v)\vert \vert u-v\vert ^{2H-2} dvdu .
\end{equation*}\end{remark}

\vskip0.5cm

The Wiener integrals $\int_{0}^{T} f(u)dZ(u)$ and $\int_{0}^{T}
g(u)dZ(u)$ are not necessarily  independent when the functions $f$
and $g$ are orthogonal in ${\cal{H}}$. A characterization of their
independence is given in the next result.

\begin{prop}
Let $f,g\in {\cal{H}}$. Then $\int_{0}^{T} f(u)dZ(u)$ and
$\int_{0}^{T} g(u)dZ(u)$ are independent if and only if
\begin{equation}
\label{indep} \langle f(\cdot ) \frac{\partial K^{H'} }{\partial u}
(\cdot, y_{1} ), g(\cdot)\frac{\partial K^{H'} }{\partial u} (\cdot, y_{2}
) \rangle _{{\cal H'}}=0 \hskip0.5cm a.e. (y_{1}, y_{2}) \in
[0,T]^{2}
\end{equation}
where ${\cal{H}'}$ is the space analogous to ${\cal{H}}$
corresponding to the Hurst parameter $H'$.
\end{prop}
{\bf Proof: } We use a result by \"Ustunel-Zakai \cite{US} (see also
Kallenberg \cite{Kal}): two multiple Wiener-It\^o integrals with respect to the
standard Wiener process $I_{n}(f)$ and $I_{m}(g)$ with $f,g$
symmetric, $f\in L^{2}[0,T]^{n}$ and $g\in L^{2}[0,T] ^{m}$ are
independent if and only if $f\otimes _{1} g=0$ a.e. on
$[0,T]^{m+n-2}$, where
\begin{equation*}
(f\otimes _{1}g) (t_{1},\ldots , t_{n-1}, s_{1}, \ldots ,s_{m-1})
=\int _{0}^{T} f(t_{1}, \ldots, t_{n-1}, t) g(s_{1}, \ldots ,
s_{n-1}, t) dt .
\end{equation*}

Let us apply the above result to
\begin{equation*}
F(y_{1}, y_{2})= 1_{[0,t]^{2}}(y_{1}, y_{2})\int _{y_{1}\vee
y_{2}}^{T} f(u)\frac{\partial K^{H'} }{\partial u} (u ,y_{1} )
\frac{\partial K^{H'} }{\partial u} (u,y_{2} ) du
\end{equation*}
and
\begin{equation*}
G(y_{1}, y_{2})= 1_{[0,t]^{2}}(y_{1}, y_{2})\int _{y_{1}\vee
y_{2}}^{T} g(u)\frac{\partial K^{H'} }{\partial u} (u ,y_{1} )
\frac{\partial K^{H'} }{\partial u} (u,y_{2} ) du.
\end{equation*}
Then
\begin{eqnarray*}
(F\otimes _{1}G) (y_{1}, y_{2}) &=& \int_{0}^{T}ds \\
&&\times \int _{y_{1}\vee s}^{T} f(u)\frac{\partial K^{H'}
}{\partial u} (u ,y_{1} ) \frac{\partial K^{H'} }{\partial u} (u,s)
du \int _{y_{2}\vee s}^{T} g(v)\frac{\partial K^{H'} }{\partial v}
(v ,y_{2} ) \frac{\partial K^{H'} }{\partial v} (v,s)dv\\
&=& c(H) \int_{0}^{T} \int _{0}^{T} f(u)g(v)\frac{\partial K^{H'}
}{\partial u} (u ,y_{1} ) \frac{\partial K^{H'} }{\partial v} (v
,y_{2} ) \vert u-v\vert ^{2H'-2}dvdu
\end{eqnarray*}
and the conclusion follows easily. \qed

\vskip0.5cm

As an immediate consequence, we obtain

\begin{corollary}
If $f\otimes g =0$ a.e. on $[0,T] ^{2}$ then the random variables
$\int_{0}^{T} f(u)dZ(u)$ and $\int_{0}^{T} g(u)dZ(u)$ are
independent.
\end{corollary}

\vskip0.5cm

\begin{remark}
The construction of Wiener integrals with respect to the Rosenblatt
process allows to consider associated Ornstein-Uhlenbeck processes.
This has been done in \cite{MaTu} for general Hermite processes of
order $k$ following the argument in \cite{CKM} for the fBm case. It
can be showed that the equation
 \begin{equation}
 \label{lange}
X_{t}= \xi  -\lambda \int_{0}^{t}X_{s}ds +\sigma Z (t), \hskip0.5cm
t\geq 0,
 \end{equation}
 where $\sigma , \lambda > 0$ and the initial condition $\xi$ is a random variable
  in  $L^{0}(\Omega)$ has an unique solution that can be represented
  as
  \begin{equation*}
X^{\xi }(t)= e^{-\lambda t} \left( \xi + \sigma \int_{0}^{t}
e^{\lambda u} dZ(u) \right), \hskip0.5cm t\geq 0.
\end{equation*}
where the stochastic integral above exists in the Wiener sense. When
the initial condition is $\xi = \sigma \int_{-\infty }^{0}
e^{\lambda u} dZ (u)$ (in \cite{CKM} the integrals are considered on the whole real line), the solution of (\ref{lange}) can be written
as
\begin{equation}
\label{fOU}X(t)= \sigma \int_{-\infty }^{t} e^{-\lambda ( t-u)}
dZ(u)
\end{equation}
and it is called the {\em stationary Rosenblatt  Ornstein-Uhlenbeck
process. }

\end{remark}

\begin{remark}
 The Non Central Limit Theorem given by \cite{DM}, \cite{Taqqu79} can
 be  extended to Wiener integrals  (see \cite{MaTu}). More precisely,
 under suitable assumptions on the deterministic function
 $f\in{\cal{H}}$ one obtains that the sequence
 \begin{equation*}
\frac{1}{n^{H}}\sum _{j\in \mathbb{Z}} f\left( \frac{j}{n}\right)
g(\xi _{j})
\end{equation*}
converges weakly when $n\to \infty$, to the Wiener integral $\int
f(u)dZ(u)$ ($g$ and $\xi_{j}$ were introduced in Section 1).
\end{remark}

\section{Infinite dimensional process and stochastic evolution
equations}

In this part we define a Hilbert-valued  Rosenblatt process and we
consider stochastic evolution equations driven by it.

Let us consider $U$ a real and separable Hilbert space and $Q$  a
nuclear, self-adjoint positive and nuclear  operator on $U$.  There
exists then a sequence $0<\lambda _{n} \searrow  0 $ of
eigenvalues of $Q$ such that $\sum _{n\geq 1} \lambda _{n} <\infty.$
Moreover the corresponding eigenvectors  form an orthonormal basis
in $U$.  We define the {\em infinite dimensional Rosenblatt process
on $U$  } as
\begin{equation}
 \label{zinf}
Z(t)=\sum _{\n \geq 0 }\sqrt { \lambda _{n}} e_{n} z_{j} (t)
 \end{equation}
 where $(z_{j})_{j\geq 0} $ is a family of real independent
 Rosenblatt  processes.

 Note that the series (\ref{zinf}) is convergent in $L^{2}(\Omega
 )$ for every $t\in [0,T]$ since
 \begin{equation*}
\mathbf{E}\left| Z(t) \right| ^{2} = \sum _{j\geq 1} \lambda _{j} \mathbf{E}(z_{j}
^{2}) = t^{2H} \sum _{j\geq 1} \lambda_{j} <\infty.
 \end{equation*}
 Note also that $Z$ {\em has covariance function $R(t,s)$ } in the
 sense that for every $u,v\in U$, and for every $s,t \in [0,T]$
 \begin{equation*}
\mathbf{E}\langle Z(t), u\rangle _{U} \langle Z(s) , v\rangle _{U}
=R(t,s)\langle Qu, v\rangle _{U}.
 \end{equation*}
This can be proved exactly  as fBm case  (see  \cite{CAT1}).
 \vskip0.5cm

 In some situations the assumption that $Q$ is nuclear is not
 convenient. For  example one cannot take  $Q$ to be the identity
 operator, that is $\lambda _{n} =1$ for every $n$. Therefore, if
 $\sum _{n} \lambda _{n} =\infty $ we will consider a bigger real
 and separable Hilbert space $U_{1} \supset U$ such  that the
 inclusion $U\subset U_{1}$ is nuclear.  Then the quantity
 \begin{equation}
\label{zinf2} Z(t)= \sum _{j} z_{j}(t) e_{j}
\end{equation}
is  well-defined stochastic  process  in $U_{1}$.

In the sequel  we will consider the infinite dimensional Rosenblatt
process to be defined by (\ref{zinf2}).

\vskip0.5cm

Following the one dimensional case, one can introduce Wiener
integrals with respect to the Hilbert-valued process $Z$.  Let $V$
be another Hilbert space Let $(\Phi _{s}, s\in [0,T]) $    a
stochastic process with valued in the space of linear operators
${\cal{L}}(U,V) $. We put for every $t\in [0,T]$
\begin{equation*}
\int _{0}^{t} \Phi _{s} dZ(s) = \sum _{j\geq 1} \int_{0}^{t} \Phi
_{s} e_{j} dZ_{j}(s)
\end{equation*}
where $ \int_{0}^{t} \Phi _{s} e_{j} dZ_{j}(s)$ is a $V$ valued
random variable. Note that the integral exists in $L^{2}(\Omega, V)$ if
\begin{equation*}
\mathbf{E} \left|  \int _{0}^{t} \Phi _{s} dZ(s) \right|_{V}  ^{2} = \sum
_{j}  \left| \Vert \Phi  e_{j} \Vert _{{\cal{H}}} \right| _{V} ^{2}
<\infty .
\end{equation*}

\begin{remark}If the integrand $ \Phi $ does not depend on time, then we find
\begin{equation*}
\mathbf{E} \left|  \int _{0}^{t} \Phi dZ(s) \right|_{V}  ^{2}= \sum _{j}
\left| \Phi e_{n} \right| _{V} ^{2} \mathbf{E} \left| \int _{0}^{t} dz_{j}(s)
\right| ^{2}=t^{2H} \sum_{j}\left| \Phi e_{n} \right| _{V} ^{2}
\end{equation*}
and it can be seen that the integral $\int _{0}^{t} \Phi dZ(s)$
exists if and only if $\Phi $ is a Hilbert-Schmidt operator.
\end{remark}

\vskip0.5cm

Now we introduce stochastic evolution equations driven by the
infinite-dimensional Rosenblatt process. Let $A : Dom (A) \subset
V\to V$  be the infinitesimal generator of the strongly continuous
semigroup $(e^{tA})_{t\in [0,T]} $. We study the equation
\begin{equation}
\label{ecu1} dX(t)= AX(t) dt + \Phi dZ(t)
\end{equation}
where $X(0)=x\in V$ and $\Phi \in  {\cal{L}}(U;V)$. We will consider
{\em mild solution } of (\ref{ecu1}), that is, (when it exists), it
can be written as
\begin{equation}
\label{mild} X(t)= e^{tA}x + \int _{0}^{t} e^{(t-s) A} \Phi dZ(s).
\end{equation}
We will not assume that $\Phi$ is Hilbert-Schmidt (although the
integral $\int \Phi dZ$ exists if and only if $\Phi$ is
Hilbert-Schmidt);  this assumption is unnecessary because, under
suitable hypothesis on $A$, the integral $\int _{0}^{t} e^{(t-s) A}
\Phi dZ(s)$ will exist even when $\Phi$ is not Hilbert-Schmidt
operator.

The method used in the fBm case will allow to prove the next
theorem.

\begin{theorem}
Let $Z$ be given by (\ref{zinf2}) with $H\in (\frac{1}{2}, 1)$.
Consider $\Phi \in  {\cal{L}}(U; V) $ and $A: Dom(A) \subset V\to V$
be a negative self-adjoint operator. Then there exists a mild
solution $X$ of the equation (\ref{ecu1}) if and only if the
operator $\Phi ^{\star } G_{H}(-A) \Phi $ is a trace class operator,
where
\begin{equation}
G_{H}(\lambda)=(\max{(\lambda,1)})^{-2H}. \label{GH}%
\end{equation}
\end{theorem}

\begin{remark}
in \cite{TTV} in the fBm case is assumed that the spectrum of $A$,
$\sigma (A) \subset -(\infty , -l ]$ with $l>0$. The situation  when
$0$ is an accumulation point of the spectrum is not treated; this
case is solved in  \cite{Dre}.
\end{remark}
{\bf Proof: } Since
\begin{equation*}
\mathbf{E}\left| \int _{0}^{t} e^{(t-s)A} \Phi dZ(s) \right| _{V} ^{2}=
c(H)\sum_{n} \int _{0}^{t} \int _{0}^{t} \langle e^{(t-u)A}\Phi
e_{n},e^{(t-v)A}\Phi e_{n}\rangle _{V} \vert u-v\vert ^{2H-2} dudv
\end{equation*}
and what it follows is a deterministic problem that can be solved as
in \cite{TTV}. \qed

\vskip0.5cm

If $\mathbb{S}^{1}$ denotes the unit circle and $A$ is the Laplacian
on the circle, we have

\begin{corollary}
Assume that  $U=V= L^{2}(\mathbb{S}^{1}) $ and $A=\Delta $ is the
Laplacian on $U$. Denote by $(e_{n}, f_{n})_{n\geq 1} $ the
eigenvectors of $\Delta $ that form an orthonormal basis in $
L^{2}(\mathbb{S}^{1}) $. Let $(q_{n})_{n}$ be a bounded sequence of
non-negative real numbers and
\begin{equation*}
Z(t)= \sum _{n} \sqrt {q_{n}} e_{n} z_{n}(t)+ \sum _{n} \sqrt
{q_{n}} f_{n} \tilde{z}_{n}(t)
\end{equation*}
with $(z_{j}, \tilde{z}_{j} )_{j}$ independent real Rosenblatt
processes.  Then (\ref{ecu1}) has an unique mild solution such that
$X(t)\in L^{2} (\Omega ,V)$ if an only if
\begin{equation*}
\sum_{n} q_{n} n^{-4H} <\infty.
\end{equation*}
\end{corollary}

\vskip0.5cm

\section{Pathwise stochastic calculus}
At this point, we will start to develop a stochastic integration
theory with respect to the Rosenblatt process. In general, for
processes that are not semimartingales, the It\^o's theory cannot be
applied. One needs generalized alternative ways to integrate
stochastically with respect to such processes. In general these
generalized method are essentially of two types: the first is the
pathwise type calculus and (here we included the rough path analysis
\cite{QL} and the stochastic calculus via regularization \cite{RV1})
and the second type is Malliavin calculus and the Skorohod
integration theory \cite{N}.  In general the pathwise type calculus
is connected to the trajectorial regularity  and/or the covariance
structure of the integrator process. The Malliavin calculus instead
is very related to the Gaussian character of the driven process.

Since the Rosenblatt process with $H>\frac{1}{2}$ has zero quadratic
variation (see \cite{RV2}) and regular paths (H\"older continuous of order $H-\varepsilon$), the pathwise calculus
can be naturally applied to construct stochastic integrals with
respect to it. Here we choose to use the approach  of Russo and
Vallois.  Let us  list first the main ingredients of the stochastic
calculus via regularization.

Let $(X_{t})_{t\geq 0}$ and and $(Y_{t}) _{t\geq 0}$ continuous
processes. We introduce, for every $t$,
\begin{equation*}
I^{-} (\varepsilon, Y, dX) =\int _{0}^{t} Y_{s}
\frac{X_{s+\varepsilon }-X_{s}} {\varepsilon } ds, \hskip0.2cm I^{+}
(\varepsilon, Y, dX) =\int _{0}^{t} Y_{s} \frac{X_{s
}-X_{(s-\varepsilon )_{+}}} {\varepsilon } ds,
\end{equation*}

\begin{equation*}
I^{0}(\varepsilon, Y, dX) =\int _{0}^{t} Y_{s}\frac{X_{s+\varepsilon
}-X_{(s-\varepsilon)_{+}}} {2\varepsilon } ds
\end{equation*}
and \begin{equation*} C_{\varepsilon }(X,Y)(t) =\int
_{0}^{t}\frac{(X_{s+\varepsilon
}-X_{(s-\varepsilon)_{+}})(Y_{s+\varepsilon
}-Y_{(s-\varepsilon)_{+}})} {\varepsilon } ds.
\end{equation*}
Then {\em the forward, backward and symmetric integrals of $Y$ with
respect to $X$}  will be given
\begin{equation*}
\int_{0}^{t}Yd^{-}X =\lim _{\varepsilon \to 0^{+} }I^{-}
(\varepsilon, Y, dX), \hskip0.3cm \int_{0}^{t}Yd^{+}X =\lim
_{\varepsilon \to 0^{+} }I^{+} (\varepsilon, Y, dX),
\end{equation*}
and \begin{equation} \int_{0}^{t}Yd^{0}X=\lim _{\varepsilon \to
0^{+} }I^{0} (\varepsilon, Y, dX)
\end{equation}
provided that the above limits exist uniformly in probability (ucp).
{\em The covariation of $X$ and $Y$}  is defined as
\begin{equation*}
[X,Y]_{t} =ucp-\lim _{\varepsilon \to 0^{+} }C_{\varepsilon
}(X,Y)(t).
\end{equation*}
If $X=Y$ we denote $[X,X]=[X]$ and when $[X]$ exists then $X$ is
said to be a finite quadratic variation process. When $[X]=0$, then
$X$ is called a zero quadratic variation process.

\vskip0.5cm

The Rosenblatt process is clearly a zero quadratic variation process
since
\begin{equation*}
\mathbf{E} C_{\varepsilon }(Z,Z)(t)= \mathbf{E} \int _{0}^{t} \frac {1}{\varepsilon }
(X_{s+\varepsilon }-X_{s}) ^{2} ds =t\varepsilon ^{2H-1} \to
_{\varepsilon \to 0} 0.
\end{equation*}
Therefore the stochastic calculus via regularization can be directly
applied to it. Precisely, it follows from Proposition 4.2 of the
Russo and Vallois survey \cite{RVsurvey} that every $f\in
C^{2}(\mathbb{R})$, the integrals $$\int_{0}^{t} f'(X) d^{-}X,
\hskip0.3cm \int_{0}^{t} f'(X) d^{+}X, \hskip0.3cm \int_{0}^{t}
f'(X) d^{0}X$$ exist and are equal  and we have the It\^o's formula
\begin{equation}
\label{ito3} f(X_{t}) =f(X_{0}) + \int _{0}^{t} f'(X)d^{0}X.
\end{equation}

\begin{remark}An immediate consequence of the existence of the quadratic
variation of the Rosenblatt process   is the existence and
uniqueness of the solution of a Stratonovich stochastic differential
equation driven by $Z$. Concretely, if $\sigma: \mathbb{R} \to
\mathbb{R}$ and $b: [0,T] \times \mathbb{R} \to \mathbb{R}$ satisfy
some regularity assumptions and $V$ is a locally bounded variation
process, then the equation
\begin{equation}
\label{stra} dX(t)= \sigma (X(t))  d^{0} Z(t) + b(t, X(t))dV(t)
\end{equation}
with $X(0)=G$ where $G$ is an arbitrary random variable, has an
unique solution (see \cite{RV2} for the definition of the solution).
\end{remark}

\vskip0.5cm

\section{Skorohod integral with respect to the Rosenblatt process}
In this part we define a divergence integral with respect to
$(Z(t))_{t\in [0,T]}$.  Constructing generalized Skorohod integrals
with respect to processes that are not necessarily Gaussian or
semimartingales constitutes a frequent  topic. For results in this
direction, we refer among others, to \cite{PriTu}, \cite{JS},
\cite{MV}, \cite{GNT} or \cite{KRT}.

We will need some basic elements of the Malliavin calculus  with
respect to a Wiener process $(W_{t})_{t\in [0,T]}$. By ${\cal{S}}$
we denote the class of smooth random variables of the form
\begin{equation}
 \label{a111}
F=f\left( W_{t_{1}}, \ldots , W_{t_{n}}\right), \hskip0.5cm t_{1},
\ldots , t_{n}\in [0,T]
\end{equation}
where  $f\in C^{\infty }_{b} (\mathbb{R}^{n})$. If $F $ is of the
form (\ref{a111}), its Malliavin derivative is defined as
\begin{equation*}
D_{t}F=\sum_{i=1}^{n}  \frac{\partial f}{\partial x_{i}}\left(
W_{t_{1}}, \ldots , W_{t_{n}}\right)1_{[0,t_{i}]}(t), \hskip0.5cm
t\in [0,T].
\end{equation*}
The operator $D$ is an unbounded closable operator and it can be
extended to the closure of ${\cal{S}}$ (denoted $\mathbb{D}^{k,p}$,
$k\geq 1$ integer, $p\geq 2$) with respect to the norm
\begin{equation*}
\Vert F\Vert ^{p}_{k,p}=\mathbf{E}\vert F\vert ^{p} + \sum
_{j=1}^{k} \mathbf{E}\Vert D^{(j)}F\Vert ^{p} _{L^{2}([0,T]^{j})},
\hskip0.5cm F\in {\cal{S}}, k\geq 1, p\geq 2
\end{equation*}
where the $j$th derivative $D^{(j)}$ is defined by iteration.

The Skorohod integral $\delta $ is the adjoint of $D$. Its domain is
\begin{equation*}
Dom (\delta) =\{ u\in L^{2}([0,T]\times \Omega )/ \left| \mathbf{E}
\int _{0}^{T}u_{s}D_{s}Fds \right| \leq C\Vert F\Vert _{2}\}
\end{equation*}
and $D$ and $\delta $ satisfy the duality relationship
\begin{equation}
\label{dua} \mathbf{E}\left( F\delta (u)\right)
=\mathbf{E}\int_{0}^{T} D_{s}Fu_{s}ds, \hskip0.5cm F\in {\cal{S}},
u\in Dom (\delta).
\end{equation}
 We
define $\mathbb{L}^{k,p}=L^{p}\left( %
\left[ 0,T\right] ;\mathbb{D}^{k,p}\right) $. Note that   $\mathbb{L}%
^{k,p}\subset Dom(\delta )$. We denote $\delta (u)= \int _{0}^{T}
u_{s} \delta W_{s}$. We will need the integration by parts formula
\begin{equation}
\label{ip1}
F\delta (u) = \delta (Fu) + \int _{0}^{T} D_{s}F u_{s}
\end{equation}
if $F\in \mathbb{D}^{1,2}$ and $u\in \mathbb{L}^{1,2}.$

\vskip0.5cm

We also mention that the Skorohod integral with respect to the fBm
$B^{H}$ with Hurst parameter $H>\frac{1}{2}$ is defined through a
transfer operator
\begin{equation}
\label{trans} \int _{0}^{T} g_{s}dB^{H}_{s}= \int _{0}^{T} \int
_{s}^{T} g_{r} \frac{\partial K^{H}}{\partial r} (r,s) dr dW_{s}
\end{equation}
where the integral in the right side above is a Skorohod integral
with respect to $W$. moreover $g$ is Skorohod integrable with
respect to $b^{H}$ if the quantity $\int _{s}^{T} g_{r}
\frac{\partial K^{H}}{\partial r} (r,s) dr $ is Skorohod integrable
with respect to $W$.

\begin{definition}
Let us consider a square integrable stochastic process $(g_{s})
_{s\in [0,T]}$.  Following (\ref{wiener1}) and (\ref{trans})  we
define its Skorohod integral with respect to $Z$ by
\begin{eqnarray}
&& \int_{0}^{T}  g_{s}dZ(s)\nonumber \\
&:=& \int _{0}^{T} \int_{0}^{T} I(g) (y_{1},
y_{2}) dB(y_{1}) dB(y_{2}) \nonumber\\
&=& \label{sko1}\int _{0}^{T} \int_{0}^{T}\left( \int _{y_{1}\vee
y_{2}}^{T} g(u)\frac{\partial K^{H'} }{\partial u} (u ,y_{1} )
\frac{\partial K^{H'} }{\partial u} (u,y_{2} ) du\right) dB(y_{1})
dB(y_{2}).
\end{eqnarray}
We will say that a process $g$ is Skorohod integrable with respect
to $Z$ if the process $Ig \in Dom \delta ^{(2)} $, where $\delta
^{(2)} $ is the double Skorohod integral with respect to the
Brownian motion $B$.
\end{definition}

We refer to \cite{NZ} for the study of double (and multiple)
Skorohod integrals.

\begin{remark}
Note that the Skorohod integral coincide with the Wiener integral if
the integrand $g$ is a deterministic function in ${\cal{H}}$.
Another Skorohod integral with respect to $Z$ has been introduced in
\cite{KRT} as the adjoint of some Malliavin derivative with respect
to $Z$ but this integral does not coincide with the Wiener integral
for deterministic integrands.
\end{remark}

Next lemma gives a condition that ensures the Skorohod
integrability.

\begin{lemma}
Let $g\in L^{2} (\Omega; {\cal{H}})$ be such that $g\in
\mathbb{L}^{2,2}$ and
\begin{equation}
\label{cond-skp} \mathbf{E}  \int _{0}^{T} \int _{0}^{T} \Vert D_{x_{1},
x_{2}}g\Vert ^{2} _{ {\cal{H}}}dx_{1}dx_{2} <\infty.
\end{equation}
Then $g $ is Skorohod integrable with respect to $Z$ and
\begin{equation}
\label{ineq1} \mathbf{E}\left|\int_{0}^{T} g_{s} \delta Z(s) \right| ^{2}
\leq cst.\left[\mathbf{E} \Vert g\Vert ^{2}_{{\cal{H}}}+\mathbf{E}  \int _{0}^{T} \int
_{0}^{T} \Vert D_{x_{1}, x_{2}}g\Vert ^{2} _{
{\cal{H}}}dx_{1}dx_{2}\right] .
\end{equation}
\end{lemma}
{\bf Proof: } We use Meyer's inequality for the double Skorohod
integral (see \cite{NZ}, pag. 320) and we obtain
\begin{eqnarray*}
\mathbf{E}\left| \int_{0}^{T} g_{s} \delta Z(s) \right| ^{2}&\leq & cst.
\left[ \mathbf{E} \int _{0}^{T} \int _{0}^{T}I(g)(y_{1}, y_{2}) ^{2}
dy_{1}dy_{2} \right.
\\
&&\left. + \mathbf{E}   \int _{0}^{T} \int _{0}^{T} \int _{0}^{T} \int
_{0}^{T}\left( D_{x_{1}, x_{2} } I(g) (y_{1}, y_{2}) \right) ^{2}
dx_{1}dx_{2} dy_{1}dy_{2} \right] \\
&=& cst. \left[ \mathbf{E}  H(2H-1)\int _{0}^{T} \int _{0}^{T} g(u)g(v)\vert
u-v\vert ^{2H-2} dvdu \right. \\
&&\left. + \int _{0}^{T} \int _{0}^{T}dx_{1}dx_{2} \left(  \int
_{0}^{T} \int _{0}^{T}D_{x_{1},x_{2}}g(u) D_{x_{1},x_{2}}g(v)\vert
u-v\vert ^{2H-2} dvdu\right)  \right] \\
&=& cst.\left[\mathbf{E} \Vert g\Vert ^{2}_{{\cal{H}}}+\mathbf{E}  \int _{0}^{T} \int
_{0}^{T} \Vert D_{x_{1}, x_{2}}g\Vert ^{2} _{
{\cal{H}}}dx_{1}dx_{2}\right].
\end{eqnarray*} \qed

\begin{corollary}
If  $g\in L^{2} (\Omega; \left| {\cal{H}}\right|)$ be such that
$g\in \mathbb{L}^{2,2}$ and
\begin{equation}
\label{cond-sko} \mathbf{E}  \int _{0}^{T} \int _{0}^{T} \Vert D_{x_{1},
x_{2}}g\Vert ^{2} _{ \left| {\cal{H}}\right| }dx_{1}dx_{2} <\infty.
\end{equation}
Then $g $ is Skorohod integrable with respect to $Z$ and
\begin{equation}
\label{ineq2} \mathbf{E}\left|\int_{0}^{T} g_{s} \delta Z(s) \right| ^{2}
\leq cst.\Vert g\Vert ^{2}
\end{equation}
where
\begin{equation*}
\Vert g\Vert ^{2}=\left[\mathbf{E} \Vert g\Vert ^{2}_{\left| {\cal{H}}\right| }+\mathbf{E} \int
_{0}^{T} \int _{0}^{T} \Vert D_{x_{1}, x_{2}}g\Vert ^{2} _{ \left|
{\cal{H}}\right| }dx_{1}dx_{2}\right].
\end{equation*}
\end{corollary}

\begin{example}
The Rosenblatt process $Z$ is Skorohod integrable with respect to
$Z$ and
\begin{equation*}
\mathbf{E}\left| \int_{0}^{T} g_{s} \delta Z(s) \right| ^{2}\leq  cst.\int
_{0}^{T} \int _{0}^{T} R(u,v) \vert u-v\vert ^{2H-2} dudv .
\end{equation*}
\end{example}
{\bf Proof: } We treat  the two terms in the right side of
(\ref{ineq1}). Clearly
\begin{equation*}
\mathbf{E}\Vert Z\Vert ^{2} _{{\cal{H}}} =cst.\int _{.0}^{T} \int _{0}^{T}
R(u,v) \vert u-v\vert ^{2H-2} dudv
\end{equation*}
We have, for every $x_{1}, x_{2} \in [0,T]$,
\begin{equation*}
D_{x_{1},x_{2}} Z(u)= 2 d(H)1_{[0,u] ^{2} }(x_{1}, x_{2} ) \int
_{x_{1} \vee x_{2}} ^{u} \frac{\partial K^{H'} }{\partial u'}
(u',x_{1})\frac{\partial K^{H'} }{\partial u'} (u',x_{2})du'
\end{equation*}
and then (note also that $ \frac{\partial K^{H'} }{\partial t}(t,s)$ is positive and therefore we omitt the absolute value at a certain point)
\begin{eqnarray*}
&& \mathbf{E}  \int _{0}^{T} \int _{0}^{T} \Vert D_{x_{1}, x_{2}}g\Vert ^{2}
_{ \left| {\cal{H}}\right| }dx_{1}dx_{2} \\
&=&
 \int _{0}^{T} \int _{0}^{T}dx_{1}dx_{2} \int _{x_{1}
\vee x_{2}}^{T} \int _{x_{1}
\vee x_{2}}^{T}\vert u-v\vert ^{2H-2} dudv\\
 &&\left| \int _{x_{1}
\vee x_{2}} ^{u} \frac{\partial K^{H'} }{\partial u'}
(u',x_{1})\frac{\partial K^{H'} }{\partial u'} (u',x_{2})du' \int _{x_{1}
\vee x_{2}} ^{v} \frac{\partial K^{H'} }{\partial v'}
(v',x_{1})\frac{\partial K^{H'} }{\partial v'} (v',x_{2})dv'\right| \\
&=& \int_{0}^{T} \int_{0}^{T}\vert u-v\vert ^{2H-2} dudv \int
_{0}^{u} \int _{0}^{v} \left(     \int _{0} ^{u'\wedge
v'}\frac{\partial K^{H'} }{\partial u'} (u',x_{1})\frac{\partial
K}{\partial v'} (v',x_{1}) dx_{1} \right) ^{2}\\
&=& cst.\int _{0}^{T} \int _{0}^{T} R(u,v) \vert u-v\vert ^{2H-2}
dudv .
\end{eqnarray*} \qed

\vskip0.5cm

We finish the section by a result on the continuity of the
indefinite Skorohod integral process. This shows that the indefinite
keeps the same order of H\"older regularity as the Rosenblatt
process.
\begin{prop}
Let $g\in \mathbb{L}^{2,p}$ such that
\begin{equation*}
\sup _{r} \Vert g_{r} \Vert _{2,p} \leq \infty.
\end{equation*}
Then the indefinite Skorohod integral process $\left(
X_{t}=\int_{0}^{t} g_{s} \delta Z(s), t\in [0,T] \right) $ admits  a
H\"older continuous version of order $\delta <H.$
\end{prop}
{\bf Proof: } We can write
\begin{eqnarray*}
X_{t}-X_{s}&=& \int_{s}^{t} \int_{s} ^{t}\left( \int _{y_{1}\vee
y_{2} }^{t} \frac{\partial K^{H'} }{\partial u} (u ,y_{1} )
\frac{\partial K^{H'} }{\partial
u} (u,y_{2} ) du\right) dB(y_{1})dB(y_{2})  \\
&& + 2 \int_{0}^{s} \int _{s}^{t} \left( \int _{y_{2} }^{t}
\frac{\partial K^{H'} }{\partial u} (u ,y_{1} )  \frac{\partial
K^{H'} }{\partial u} (u,y_{2} ) du\right) dB(y_{1} ) dB(y_{2}) \\
&&+\int_{0}^{s} \int _{0}^{s} \left( \int _{s }^{t} \frac{\partial
K^{H'} }{\partial u} (u ,y_{1} )  \frac{\partial K^{H'} }{\partial
u} (u,y_{2} ) du\right) dB(y_{1} ) dB(y_{2})\\
&:=& J_{1} +2J_{2} + J_{3}.
\end{eqnarray*}
Then
\begin{equation*}
\mathbf{E}\left| X_{t} -X_{s} \right| ^{p} \leq c(p) \mathbf{E}\left( J_{1}^{p} +
J_{2} ^{p} + J_{3} ^{p}\right).
\end{equation*}
By Meyer's inequality (\cite{NZ}, pag. 320)
\begin{eqnarray*}
\mathbf{E}\vert J_{1} \vert ^{p} &\leq & c(p)  \left| \int_{s}^{t} \int
_{s}^{t} \mathbf{E} \left(  \int _{y_{1}\vee y_{2} }^{t} g_{u}\frac{\partial
K^{H'} }{\partial u} (u ,y_{1} ) \frac{\partial K^{H'} }{\partial u}
(u,y_{2} ) du  \right)^{2} dy_{1} dy_{2} \right| ^{\frac{p}{2}} \\
&&+c(p)\mathbf{E} \left| \int _{s}^{t} \int _{s}^{t} \int _{0}^{t} \int
_{0}^{t}\left[    D _{x_{1}, x_{2}}\int _{y_{1}\vee y_{2}
}^{t} g_ {u}\frac{\partial K^{H'} }{\partial u} (u ,y_{1} )
\frac{\partial K^{H'} }{\partial u} (u,y_{2} ) du  \right]^{2}
dx_{1}dx_{2} dy_{1}
dy_{2}\right| ^{\frac{p}{2}} \\
&=&c(p,H)\left| \mathbf{E}\int_{s}^{t} \int _{s}^{t} \vert g(u)g(v)\vert
\vert u-v\vert ^{2H-2}\right| ^{\frac{p}{2}} \\
&&+c(p,H)\mathbf{E}\left| \int _{0}^{t} \int _{0}^{t}dx_{1} dx_{2}
\int_{s}^{t} \int_{s}^{t} \vert D _{x_{1}, x_{2}}g_ {u}D
_{x_{1}, x_{2}} g_ {v} \vert \vert u-v\vert ^{2H-2} dvdu \right|
^{\frac{p}{2}}\\
&\leq & c(p,H)\sup _{r} \Vert g_{r} \Vert _{2,p}^{p} \left|
\int_{s}^{t}\int_{s}^{t}\vert u-v\vert ^{2H-2} dvdu\right|
^{\frac{p}{2}} \\
&=& c(p,H)\sup _{r} \Vert g_{r} \Vert _{2,p}^{p}(t-s)^{pH} .
\end{eqnarray*}
In a similar way, we can find the same bound for the terms $J_{2}$
and $J_{3}$ (see also \cite{AMN}, proof of Proposition 1). The conclusion will following by the Kolmogorov's
continuity criterium. \qed

\vskip0.5cm

\section{The relation between the pathwise and the Skorohod
integrals} Let $g$ a stochastic process. Recall that its forward
integral with respect to $Z$ is the limit ucp as $\varepsilon \to 0$
of
\begin{equation}\label{a1}
I^{-}(\varepsilon, g, dZ)= \frac{1}{\varepsilon}  \int _{0}^{T}
g_{s} (Z(s+\varepsilon ) -Z(s) ) ds =\frac{1}{\varepsilon}  \int
_{0}^{T} g_{s}  \delta ^{(2)} \left( f_{s+\varepsilon }(\cdot, \ast
) -f_{s}(\cdot, \ast ) \right) ds
\end{equation}
where the kernel $f_{s}$ is given by
\begin{equation}
\label{fs} f_{s}(x,y)= d(H)1_{[0, s] ^{2}}  (x,y) \int _{x\vee y}
^{s} \frac{\partial K^{H'} } {\partial u } (u,x) \frac{\partial K^{H'} }
{\partial u } (u,y)du.
\end{equation}

\vskip0.5cm

We will need a formula by \cite{NZ}: if $F\in \mathbb{D} ^{2,2} $ ,
$u\in L^{2}([0,T]^{2}\times \Omega ) $ such that for every $s$,
$u(\cdot , s) \in Dom (\delta )$, then $Fu \in Dom (\delta ^{(2)} )
$ and
\begin{equation}
\label{ip2} F\delta ^{(2)}(u) = \delta ^{(2)} (Fu) + 2 \int _{0}^{T}
D_{\alpha } F \delta (u(\cdot , \alpha ) ) d\alpha -\int_{0}^{T}
\int_{0}^{T} D_{\alpha , \beta } ^{(2)} F u(\alpha , \beta ) d\alpha
d\beta .
\end{equation}

\vskip0.5cm

We apply relation (\ref{ip2}) to (\ref{a1}) and we obtain
\begin{eqnarray}
I^{-}(\varepsilon, g, dZ)&=&\frac{1}{\varepsilon} \int _{0}^{T}
\delta ^{(2)} \left(  g_{s}\left( f_{s+\varepsilon }(\cdot, \ast )
-f_{s}(\cdot, \ast ) \right) \right) ds \nonumber  \\
&&+\frac{2}{\varepsilon} \int _{0}^{T}\int_{0}^{T} D_{\alpha }g_{s}
\delta \left( f_{s+\varepsilon }(\cdot, \alpha )
-f_{s}(\cdot, \alpha ) \right)d\alpha ds \nonumber \\
&&-\frac{1}{\varepsilon} \int _{0}^{T}\int _{0}^{T}\int _{0}^{T}
D_{\alpha , \beta } ^{(2)}g_{s}\left( f_{s+\varepsilon }(\beta,
\alpha ) -f_{s}(\beta , \alpha ) \right) d\beta d\alpha ds.
\label{nou1}
\end{eqnarray}
We can already observe, besides the first divergence type term, the
appearance  of two trace terms. Recall that in the fBm case the
corresponding term $I^{-}(\varepsilon, g, dB^{H})$ can be decomposed
 in a divergence term plus a only  a trace term.

\begin{definition}
We say that a stochastic process $g\in \mathbb{L}^{1,2}$ {\em admits
a trace of order 1 }if
\begin{equation}
\label{tr1}\frac{1}{\varepsilon} \int _{0}^{T}\int_{0}^{T} D_{\alpha
}g_{s} \delta \left( f_{s+\varepsilon }(\cdot, \alpha )
-f_{s}(\cdot, \alpha ) \right)d\alpha ds
\end{equation}
converges in probability as $\varepsilon \to 0$. The limit will be
denoted by $Tr^{(1)}(D^{(1)}g)$.

\vskip0.5cm

We say that a stochastic process $g\in \mathbb{L}^{2,2}$ {\em admits
a trace of order 2 } if \begin{equation} \label{tr2}
\frac{1}{\varepsilon} \int _{0}^{T}\int _{0}^{T}\int _{0}^{T}
D_{\alpha , \beta } ^{(2)}g_{s}\left( f_{s+\varepsilon }(\beta,
\alpha ) -f_{s}(\beta , \alpha ) \right) d\beta d\alpha ds
\end{equation}
converges in probability as $\varepsilon \to 0$. The limit will be
denoted by $Tr^{(2)}(D^{(2)}g)$. \label{trace}
\end{definition}

We have the following relation between the divergence and the
pathwise integral.

\begin{theorem}\label{SS}
Let $g\in \mathbb{L} ^{2,2}$ such that \begin{equation*}
 \mathbf{E} \Vert
g\Vert ^{2}_{\left| {\cal{H}}\right| }+\mathbf{E} \int _{0}^{T} \int _{0}^{T}
\Vert D_{x_{1}, x_{2}}g\Vert ^{2} _{ \left| {\cal{H}}\right|
}dx_{1}dx_{2}<\infty.
\end{equation*}
Assume that $g$ has traces of order 1 and 2. Then $g$ is forward
integrable with respect to $Z$ and it holds
 \begin{equation}
\label{sko-path} \int_{0}^{T} g_{s} d^{-} Z(s) = \int _{0}^{T} g_{s}
\delta Z(s) +2 Tr ^{(1)} (D^{(1)}g) - Tr^{(2)}(D^{(2) }g).
\end{equation}\end{theorem}
{\bf Proof: } By (\ref{nou1}) and Definition \ref{trace}, it
suffices to show that the term
\begin{equation*}
A_{\epsilon}= \frac{1}{\varepsilon} \int _{0}^{T} \delta ^{(2)}
\left( g_{s}\left( f_{s+\varepsilon }(\cdot, \ast ) -f_{s}(\cdot,
\ast ) \right) \right) ds
\end{equation*}
converges to
\begin{equation*}
\int_{0}^{T} g_{s} \delta Z(s)
\end{equation*}
in $L^{2}(\Omega )$ as $\varepsilon \to 0$.  \vskip0.5cm

We can write, by Fubini,
\begin{eqnarray*}
A_{\varepsilon}&=& \frac{1}{\varepsilon} \int_{0}^{T} ds
\int_{0}^{T} \int _{0}^{T} g_{s} \left( f_{s+\varepsilon }(y_{1},
y_{2}) -f_{s }(y_{1}, y_{2})\right)  dB(y_{1})dB(y_{2}) \\
&=& \int_{0}^{T} \int_{0}^{T}dB(y_{1})dB(y_{2}) \int_{y_{1}\vee
y_{2}} ^{T} g^{\varepsilon } (u) \frac{\partial K^{H'} }{\partial u} (u,
y_{1})\frac{\partial K^{H'} }{\partial u} (u, y_{2})du \\
&=&\int_{0}^{T}\int_{0}^{T}I(g^{\varepsilon})(y_{1}, y_{2})
dB(y_{1})dB(y_{2}) = \int _{0}^{T} g^{\varepsilon } _{s} \delta Z(s)
\end{eqnarray*}
where we denoted by
\begin{equation}
\label{geps} g^{\varepsilon }(u)=\frac{1}{\varepsilon} \int
_{u-\varepsilon}^{u} g_{s}ds.
\end{equation}
By using (\ref{ineq2}), it is sufficient to check that
\begin{equation*}
g^{\varepsilon } \to _{\varepsilon \to 0} g \mbox{ in } L^{2}(\Omega
; {\cal{H}})
\end{equation*}
and \begin{equation*} \int_{0}^{T}\int_{0}^{T} \mathbf{E}\Vert D
_{x_{1}, x_{2} } (g^{\varepsilon } -g)\Vert _{{\cal{H}}}^{2}
dx_{1}dx_{2} \to _{\varepsilon \to 0}0.
\end{equation*}
We will show that
\begin{equation}
\label{a2} \Vert g^{\varepsilon } \Vert _{\left| {\cal{H}}\right| }
\leq c(H)\Vert g \Vert _{\left| {\cal{H}}\right| }
\end{equation}
and
\begin{equation}
\label{a3}  \int _{0}^{T} \int _{0}^{T} \Vert D_{x_{1},
x_{2}}g^{\varepsilon} \Vert ^{2} _{ \left| {\cal{H}}\right|
}dx_{1}dx_{2}\leq c(H) \int _{0}^{T} \int _{0}^{T} \Vert D_{x_{1},
x_{2}}g\Vert ^{2} _{ \left| {\cal{H}}\right| }dx_{1}dx_{2}.
\end{equation}
The bound (\ref{a2}) has been proved in \cite{AN}, proof of
Proposition 3, Step 1. Concerning the bound (\ref{a3}), we can write
\begin{eqnarray*}
&& \int _{0}^{T} \int _{0}^{T} \Vert D_{x_{1},
x_{2}}g^{\varepsilon}
\Vert ^{2} _{ \left| {\cal{H}}\right| }dx_{1}dx_{2}\\
&= & c(H) \int _{0}^{T} \int _{0}^{T}dx_{1}dx_{2}\int _{0}^{T} \int
_{0}^{T}\left| D _{x_{1}, x_{2}}g^{\varepsilon}_{u} \right| \left| D _{x_{1}, x_{2}}g^{\varepsilon}_{v}\right|
 \vert u-v\vert ^{2H-2}dudv\\
&\leq & c(H) \frac{1}{\varepsilon ^{2}} \int _{0}^{T} \int
_{0}^{T}dx_{1}dx_{2}\int _{0}^{T} \int _{0}^{T} dudv \vert u-v\vert
^{2H-2}\int_{v-\varepsilon }^{v}\int_{u-\varepsilon }^{u}ds ds' \left| D _{x_{1}, x_{2}}g_{s}D _{x_{1}, x_{2}}g_{s'}\right| \\
\\
&\leq & c(H)\int _{0}^{T} \int _{0}^{T}dx_{1}dx_{2}\int _{0}^{T}
\int _{0}^{T}dsds' \left| D _{x_{1}, x_{2}}g_{s}D _{x_{1}, x_{2}}g_{s'}\right| \left( \frac{1}{\varepsilon ^{2}}
\int_{s}^{s+\varepsilon } \int_{s'}^{s'+\varepsilon } \vert u-v\vert
^{2H-2}dudv \right) .
\end{eqnarray*}
It follows from \cite{AN}, proof of Proposition 3, Step 1, that
\begin{equation*}
\frac{1}{\varepsilon ^{2}} \int_{s}^{s+\varepsilon }
\int_{s'}^{s'+\varepsilon } \vert u-v\vert ^{2H-2}dudv \leq c(H)
\vert s-s'\vert ^{2H-2}
\end{equation*}
and thus (\ref{a3}) follows.

Now we can finish the proof proceeding as in \cite{AN}, proof of
Proposition 3, Step 3. Consider a sequence $g^{n}$ of simple
processes of the form $g^{n} =\sum _{i=0}^{n-1} F_{i} 1_{(t_{i},
t_{i+1}] }$ with $F_{i} \in {\cal{S}}$ and $t_{i}\in [0, T]$ such
that $\Vert g_{n} -g \Vert \to 0$ in $L^{2}(\Omega )$ when $n\to
\infty$ (the existence of a  such sequence follows easily by the densite of ${\cal{E}}$ in ${\cal{H}}$).  Then by (\ref{ineq2}) we have that
\begin{equation*}
\int _{0}^{T} g^{n}_{s} \delta Z(s) \to _{n\to \infty }\int _{0}^{T}
g_{s} \delta Z(s)
\end{equation*}
Denote by $g^{n, \varepsilon }$ the approximation process of the
form (\ref{geps}) associated to $g^{n}$. We can write, for any
$\varepsilon >0$ and $n\geq 1$,
\begin{eqnarray*}
&&\mathbf{E} \left| \int _{0}^{T} g^{\varepsilon } _{s} \delta Z(s) -   \int
_{0}^{T} g _{s} \delta Z(s)\right| ^{2} \leq 3\left( \mathbf{E}\left|\int
_{0}^{T} g^{\varepsilon } _{s} \delta Z(s) - \int _{0}^{T} g^{n,
\varepsilon } _{s} \delta Z(s)\right| ^{2} \right. \\
&& \left. + \mathbf{E}\left|\int _{0}^{T} g^{n,\varepsilon } _{s} \delta Z(s)
- \int _{0}^{T} g^{n} _{s} \delta Z(s)\right| ^{2}  + \mathbf{E}\left| \int
_{0}^{T} g^{n } _{s} \delta Z(s) - \int _{0}^{T} g _{s} \delta
Z(s)\right| ^{2} \right).
\end{eqnarray*}
By (\ref{a2}) and (\ref{a3}) it follows that for $n$ large enough
and for any $\delta >0$
\begin{equation*}
\mathbf{E} \left| \int _{0}^{T} g^{\varepsilon } _{s} \delta Z(s) -   \int
_{0}^{T} g _{s} \delta Z(s)\right| ^{2} \leq 3\left( \mathbf{E}\left|\int
_{0}^{T} g^{n,\varepsilon } _{s} \delta Z(s) - \int _{0}^{T} g^{n}
_{s} \delta Z(s)\right| ^{2}  +\delta \right)
\end{equation*}
and we can conclude by taking $\varepsilon \to 0$.
 \qed

\vskip0.5cm

\section{On the It\^o formula in the Skorohod sense}

We study It\^o's formula  for the Rosenblatt process in the
divergence sense. As we mentioned before, the Gaussian nature of the
integrator process is essential in the framework of the divergence
calculus and this fact can be entirely observed here. We are
actually able to prove Skorohod It\^o's formula only in two
particular cases; but more relevant than these formulas, which are
not easily tractable, is the fact that one can observe from the
computations contained here that the standard method to obtain
divergence type change of variables formulas (see e.g. \cite{N})
does not work here, in the sense that one cannot hope to obtain
It\^o's formulas that stop at $f''$.

\vskip0.5cm

 We will deduce the Skorohod It\^o formula by using
the pathwise It\^o formula. Recall that for any function $f\in
C^{2}(\mathbb{R})$
\begin{eqnarray*}
f(Z(t))&=& f(0)+ \int _{0}^{t} f'(Z(s))d^{-}Z(s) \\
&=& f(0)+ \int _{0}^{t} f'(Z(s)) \delta Z(s) + 2Tr ^{(1)} (D^{(1)}
f'(Z(s))) - Tr  ^{(2)} (D^{(2)}  f''(Z(s)))
\end{eqnarray*}
provided that the above terms exist.

\subsection{The trace of order $1$}

 Recall  that
\begin{equation*}
Tr ^{(1)} (D^{(1)} f'(Z(s))) =ucp-\lim _{\varepsilon \to 0 }
B_{\varepsilon }
\end{equation*}
where
\begin{eqnarray*}
B_{\varepsilon }&=&\frac{1}{\varepsilon } \int _{0}^{t}  ds \int
_{0}^{t} d\alpha D_{\alpha } f'(Z(s)) \delta \left( f_{s+\varepsilon
}(\cdot , \alpha ) -s_{s}(\cdot , \alpha ) \right) \\
&=& \frac{1}{\varepsilon } \int _{0}^{t}  ds \int _{0}^{t} d\alpha
f'' (Z(s)) D_{\alpha } Z(s)\delta \left( f_{s+\varepsilon }(\cdot ,
\alpha ) -f_{s}(\cdot , \alpha ) \right).
\end{eqnarray*}
The Malliavin derivative of $Z(s)$ is given by
\begin{equation}
  D_{\alpha }Z(s)=2d(H) 1_{[0,s]}(\alpha )
\left( \int_{0}^{s} \left( \int _{\alpha \vee y_{1}} ^{s} \frac{
\partial K^{H'}} {\partial u} (u, \alpha ) \frac{\partial K^{H'} } {\partial u}
(u,y_{1}) du \right) dB(y_{1}) \right) \label{dz}
\end{equation}
Thus
\begin{eqnarray*}
B_{\varepsilon }&=& \frac{2}{\varepsilon } d(H)\int _{0}^{T} ds
f''(Z(s)) \int _{0}^{s} \delta\left( f_{s+\varepsilon }(\cdot
, \alpha ) -f_{s}(\cdot , \alpha ) \right)d\alpha \\
&&  \int _{0} ^{s}\left( \int _{\alpha \vee y_{1}} ^{s} \frac{
\partial K^{H'}} {\partial u} (u, \alpha ) \frac{\partial K^{H'} } {\partial u}
(u,y_{1}) du \right) dB(y_{1})
\end{eqnarray*}
where $f_{s}$ is given by (\ref{fs}).  By using the integration by
parts formula  (\ref{ip1}) it holds that
\begin{eqnarray*}
B_{\varepsilon }&=& \frac{2}{\varepsilon }d(H)\int _{0}^{T} ds
f''(Z(s))\int _{0}^{s} d\alpha \\
&&\int _{0}^{s}\left[ \delta\left( f_{s+\varepsilon }(\cdot , \alpha
) -f_{s}(\cdot , \alpha ) \right)\int _{\alpha \vee y_{1}} ^{s}
\frac{
\partial K^{H'}} {\partial u} (u, \alpha ) \frac{\partial K^{H'} } {\partial u}
(u,y_{1}) du \right] dB(y_{1})\\
&&+ \frac{2}{\varepsilon }d(H)\int _{0}^{T} ds
f''(Z(s))\int _{0}^{s} d\alpha \\
&&\int _{0}^{s}\left[\int _{\alpha \vee y_{1}} ^{s} \frac{
\partial K^{H'}} {\partial u} (u, \alpha ) \frac{\partial K^{H'} } {\partial u}
(u,y_{1}) du (f_{s+ \varepsilon } (y_{1}, \alpha ) -f_{s} (y_{1},
\alpha ) \right] dy_{1}\\
&:= & B^{1} _{\varepsilon } + B^{2} _{\varepsilon }.
\end{eqnarray*}
We regard  first the term $B^{2} _{\varepsilon }$  because it can be
treated in the same manner for any function $f$. We can write
\begin{eqnarray*}
B^{2}_{\varepsilon }&=& \frac{2}{\varepsilon }d(H)^{2}\int _{0}^{T}
ds
f''(Z(s))\int _{0}^{s} d\alpha \\
&& \int _{0}^{s} dy_{1}\int _{\alpha \vee y_{1}} ^{s} \frac{
\partial K^{H'}} {\partial u} (u, \alpha ) \frac{\partial K^{H'} } {\partial u}
(u,y_{1}) du \\
&& \left[ 1_{[0, s+ \varepsilon ] ^{2} }(y_{1}, \alpha )\int
_{\alpha \vee y_{1}} ^{s+ \varepsilon } \frac{
\partial K^{H'}} {\partial v} (v, \alpha ) \frac{\partial K^{H'} } {\partial v}
(v,y_{1}) dv\right. \\
&&\left. -1_{[0, s ] ^{2} }(y_{1}, \alpha )\int _{\alpha \vee y_{1}}
^{s} \frac{
\partial K^{H'}} {\partial v} (v, \alpha ) \frac{\partial K^{H'} } {\partial v}
(v,y_{1}) dv\right] \\
&=& \frac{2}{\varepsilon }d(H)^{2}\int _{0}^{T} ds f''(Z(s)) \int
_{0}^{s}du \int _{s} ^{s+ \varepsilon }dv   \left( \int
_{0}^{u\wedge v} \frac{
\partial K^{H'}} {\partial u} (u, \alpha )  \frac{
\partial K^{H'}} {\partial v} (v, \alpha )d\alpha \right) ^{2}\\
&=& 2A(H) ^{2} \int _{0}^{t} dv \int _{0} ^{v}  du \vert u-v\vert
^{2H-2} \frac{1}{\varepsilon } \int _{(v-\varepsilon ) \vee u
}^{v} f''(Z(s))  ds
\end{eqnarray*}
with $A(H)= H'(2H'-1) d(H) $ and we have
\begin{eqnarray}
B^{2}_{\varepsilon }&=&  2A(H) ^{2} \int _{0}^{t} dv \int _{0} ^{v}  du \vert u-v\vert
^{2H-2}\left(   \frac{1}{\varepsilon } \int _{(v-\varepsilon )
}^{v} f''(Z(s))  ds\right) \nonumber\\
&&+ 2A(H) ^{2} \int _{0}^{t} dv \int _{v-\varepsilon } ^{v}  du \vert u-v\vert
^{2H-2}\left(   \frac{1}{\varepsilon } \int _{u }^{v} f''(Z(s))  ds\right) . \label{ita}
\end{eqnarray}

Therefore we have the convergence in $L^{1}(\Omega)$ as $\varepsilon \to 0$
\begin{equation} \label{b2}
B^{2}_{\varepsilon } \to 2A(H)^{2} \int _{0} ^{T} \int _{0}^{u}
f''(Z(v)) \vert u-v\vert ^{2H-2} dvdu=\frac{A(H) ^{2}}{2H-1} \int
_{0}^{T} \int _{0} ^{T}f''(Z(u)) u^{2H-1} du
\end{equation}
since the first summand in (\ref{ita}) converges to the limit and the second one goes to zero by the dominated convergence theorem.

The study of the term $B^{1}$ is rather difficult to be done in
general. We will restrict ourselves to its study in some particular
cases.

\vskip0.5cm {\bf The case $f(x)=x^{2}$ : }

 We have
 \begin{eqnarray*}
 B^{1}_{\varepsilon } &=& \frac{4}{\varepsilon }d(H) \int _{0}^{t} ds
 \int_{0}^{s} d\alpha \int _{0} ^{t}   \int _{0} ^{t} dB(y_{1})
 dB(y_{2}) \\
 &&\int _{\alpha \vee y_{1}} ^{s} \frac{
\partial K^{H'}} {\partial u} (u, \alpha ) \frac{\partial K^{H'} } {\partial u}
(u,y_{1}) du  (f_{s+ \varepsilon}(y_{2}, \alpha ) -f_{s} ( y_{2},
\alpha ))
\end{eqnarray*}
and by Fubini  we get
\begin{eqnarray*}
 B^{1}_{\varepsilon } &=& 4d(H)^{2}  \int _{0} ^{t} \int
 _{0} ^{t} dB(y_{1})
 dB(y_{2}) \int _{y_{2} }^{t}dv  \int_{y_{1}} ^{v}du \frac{\partial K^{H'} } {\partial u}
(u,y_{1}) \frac{\partial K^{H'} } {\partial u} (u,y_{2})\\
&&\left(  \frac{1}{\varepsilon }\int _{(v-\varepsilon) \vee   u}
^{v} ds \right) \left(  \int _{0} ^{u\wedge v} \frac{
\partial K^{H'}} {\partial u} (u, \alpha ) \frac{
\partial K^{H'}} {\partial v} (v, \alpha )d\alpha \right) \\
&=& 4d(H)^{2} H'(2H'-1)  \int _{0} ^{t} \int
 _{0} ^{t} dB(y_{1})
 dB(y_{2})\\
 &&  \int _{y_{2} }^{t}dv  \int_{y_{1}} ^{v}du  \frac{\partial K^{H'} } {\partial u}
(u,y_{1}) \frac{\partial K^{H'} } {\partial u} (u,y_{2}) \vert u-v\vert
^{2H'-2}\left(  \frac{1}{\varepsilon }\int _{(v-\varepsilon) \vee u}
^{v} ds\right)
\end{eqnarray*}
and then one can prove that
\begin{eqnarray}
&& B^{1}_{\varepsilon }\to   4d(H)^{2} H'(2H'-1)  \int _{0}
^{t} \int
 _{0} ^{t} dB(y_{1})
 dB(y_{2})\nonumber \\
&&\label{b1} \int _{y_{2} }^{t}dv  \int_{y_{1}} ^{v}du \frac{\partial K^{H'} } {\partial u}
(u,y_{1}) \frac{\partial K^{H'} } {\partial u} (u,y_{2}) \vert u-v\vert
^{2H'-2}.
\end{eqnarray}

\vskip0.5cm

{\bf The case $f(x)=x ^{3}$: } It holds by calculating first the integral $d\alpha $

\begin{eqnarray*}
B_{\varepsilon }^{1} &=& \frac{12}{\varepsilon } d(H)^{2} H'(2H'-1)
\int _{0} ^{t}
Z(s)ds \\
&& \left[ \int_{0} ^{s} \int _{0} ^{s+\varepsilon } dB(y_{1})
 dB(y_{2}) \int _{y_{1}}^{s}du \int _{y_{2}}^{s+\varepsilon }dv\frac{\partial K^{H'} } {\partial u}
(u,y_{1}) \frac{\partial K^{H'} } {\partial v} (v,y_{2})\vert u-v\vert
^{2H'-2}\right. \\
&&  \left. -\int_{0} ^{s} \int _{0} ^{s } dB(y_{1})
 dB(y_{2}) \int _{y_{1}}^{s}du \int _{y_{2}}^{s }dv\frac{\partial K^{H'} } {\partial u}
(u,y_{1}) \frac{\partial K^{H'} } {\partial v} (v,y_{2})\vert u-v\vert
^{2H'-2}\right]
\end{eqnarray*}
and the integration by parts formula for the double Skorohod
integral gives
\begin{equation*}
B_{\varepsilon }^{1}= B _{\varepsilon }^{1,1} +B _{\varepsilon
}^{1,2}+ B _{\varepsilon }^{1,3}
\end{equation*}
where
\begin{eqnarray*}
B_{\varepsilon }^{1,1} &=& \frac{12}{\varepsilon } d(H)^{2}
H'(2H'-1) \int _{0} ^{t}
ds \\
&& \left[ \int_{0} ^{s} \int _{0} ^{s+\varepsilon } dB(y_{1})
 dB(y_{2}) Z(s)\int _{y_{1}}^{s}du \int _{y_{2}}^{s+\varepsilon }dv\frac{\partial K^{H'} } {\partial u}
(u,y_{1}) \frac{\partial K^{H'} } {\partial v} (v,y_{2})\vert u-v\vert
^{2H'-2}\right. \\
&&  \left. -\int_{0} ^{s} \int _{0} ^{s } dB(y_{1})
 dB(y_{2})Z(s) \int _{y_{1}}^{s}du \int _{y_{2}}^{s }dv\frac{\partial K^{H'} } {\partial u}
(u,y_{1}) \frac{\partial K^{H'} } {\partial v} (v,y_{2})\vert u-v\vert
^{2H'-2}\right],
\end{eqnarray*}

\begin{eqnarray*}
B _{\varepsilon }^{1,2}&=& -\frac{48}{\varepsilon }d(H)^{3} H'(2H'-1) \int _{0} ^{t} ds\int _{0}^{s} d\alpha \\
&&\int_{0}^{s} \left( \int _{\alpha  \vee y_{1} } ^{s}du'
\frac{\partial K^{H'} } {\partial u'} (u', \alpha ) \frac{\partial K^{H'} }
{\partial u'}(u',
y_{1})\right)  dB(y_{1}) \\
&& \left[ \int _{0} ^{s+ \varepsilon } dB(y_{2}) \int _{\alpha }^{s}du
\int _{y_{2} }^{s+ \varepsilon  }dv\frac{\partial K^{H'} }{\partial u }(u,
\alpha ) \frac{\partial K^{H'} }{\partial v }(v, y_{2}  )\vert u-v\vert
^{2H'-2}
 \right.  \\
&& \left. -\int _{0} ^{s } dB(y_{2}) \int _{\alpha }^{s}du \int
_{y_{2} }^{s }dv\frac{\partial K^{H'} }{\partial u }(u, \alpha )
\frac{\partial K^{H'} }{\partial v }(v, y_{2}  )\vert u-v\vert ^{2H'-2}
\right]
\end{eqnarray*}
and
\begin{eqnarray*}
B^{1,3}_{\varepsilon }&=&\frac{24}{\varepsilon }d(H)^{3} H'(2H'-1)
\int _{0} ^{t}
ds  \\
&& \left[ \int  _{0} ^{s} \int _{0} ^{s+ \varepsilon } dy_{1} dy_{2}
\int _{y_{1} \vee y_{2} }^{s}  \frac{\partial K^{H'} } {\partial u'} (u',
y_{1}) \frac{\partial K^{H'} } {\partial u'} (u', y_{2}) du' \right. \\
&&\left. \int _{y_{1}} ^{s} du\int_{y_{2}} ^{s+ \varepsilon }dv
\frac{\partial K^{H'} } {\partial u} (u, y_{1}) \frac{\partial K^{H'} }
{\partial v} (v,y_{2}) \vert u-v\vert ^{ 2H'-2 }\right.\\
&&\left. -\int  _{0} ^{s} \int _{0} ^{s } dy_{1} dy_{2} \int _{y_{1}
\vee y_{2} }^{s}  \frac{\partial K^{H'} } {\partial u'} (u',
y_{1}) \frac{\partial K^{H'} } {\partial u'} (u', y_{2}) du' \right. \\
&&\left. \int _{y_{1}} ^{s} du\int_{y_{2}} ^{s }dv \frac{\partial K^{H'} }
{\partial u} (u, y_{1}) \frac{\partial K^{H'} } {\partial v} (v,y_{2})
\vert u-v\vert ^{ 2H'-2 }\right]
\end{eqnarray*}
Now,
\begin{eqnarray*}
B^{1,1} _{\varepsilon } &=& 12d(H)^{2} H'(2H'-1) \int _{0}^{t} \int
_{0} ^{t} dB(y_{1})
 dB(y_{2}) \\
 && \int _{y_{2}} ^{t}dv\int _{y_{1}}^{v} du\frac{\partial K^{H'} }
{\partial u} (u, y_{1}) \frac{\partial K^{H'} } {\partial v} (v,y_{2})
\vert u-v\vert ^{ 2H'-2 } \left( \frac{1} {\varepsilon } \int
_{(v-\varepsilon )\vee u} ^{v} Z(s) ds\right)
\end{eqnarray*}
and we have the convergence as in the proof of Theorem \ref{SS}
\begin{eqnarray}
&& B^{1,1} _{\varepsilon }\to 12d(H)^{2} H'(2H'-1)  \int
_{0}^{t} \int _{0} ^{t} dB(y_{1})
 dB(y_{2})\nonumber \\
&& \label{b11}\int _{y_{2}} ^{t}dv \int _{y_{1}}^{v} duZ_{v}\frac{\partial K^{H'} }
{\partial u} (u, y_{1}) \frac{\partial K^{H'} } {\partial v} (v,y_{2})
\vert u-v\vert ^{ 2H'-2 }
\end{eqnarray}
To treat the term $B^{1,2} _{\varepsilon }$ one needs to use again
the integration by parts formula (\ref{ip1}) and one obtains
\begin{eqnarray*}
B^{1,2}_{\varepsilon } &\to & -48 d(H)^{3}( H'(2H'-1))^{2}  \int _{0}^{t}
\int _{0} ^{t} dB(y_{1})
 dB(y_{2}) \\
 &&\int _{y_{2}} ^{t} dv\int _{y_{1}}^{v}du' \int _{y_{2}}^{v}du
\vert u-u'\vert ^{2H'-2} \vert u-v\vert ^{2H'-2} \frac{\partial K^{H'} } {\partial u'} (u',
y_{1}) \frac{\partial K^{H'} } {\partial v} (v,y_{2}) \\
&& -48d(H)^{3}( H'(2H'-1))^{3} \int_{0}^{t} dv \int _{0} ^{v} \int_{0}^{v} du'du
\vert u-v\vert ^{2H'-2}\vert u'-v\vert ^{2H'-2}\vert u-u'\vert
^{2H'-2}.
\end{eqnarray*}
Concerning $B^{1,3}_{\varepsilon }$ we similarly have
\begin{equation}
\label{b13} B^{1,3}_{\varepsilon }\to 24d(H)^{3}(H'(2H'-1))^{3} \int_{0}^{t} dv
\int _{0} ^{v} \int_{0}^{v} du'du \vert u-v\vert ^{2H'-2}\vert
u'-v\vert ^{2H'-2}\vert u-u'\vert ^{2H'-2}.
\end{equation}

\vskip0.5cm

\subsection{The trace  of order 2 }

Recall that
\begin{equation*}
Tr^{(2)} \left(    D^{(2)} f'(Z(s)) \right) =ucp-\lim _{\varepsilon
\to 0} C_{\varepsilon }
\end{equation*}
where
\begin{eqnarray*}
C_{\varepsilon }   &=&\frac{1} {\varepsilon } \int _{0} ^{t} ds
\int_{0}^{t} \int_ {0}^{t} D_{\alpha , \beta } f'(Z(s)) \left(
f_{s+\varepsilon }(\alpha , \beta ) -f_{s }(\alpha , \beta )\right)
d\alpha d\beta \\
&=& \frac{1} {\varepsilon } \int _{0} ^{t} ds\int_{0}^{t}
\int_{0}^{t}\left( f_{s+\varepsilon }(\alpha , \beta ) -f_{s
}(\alpha , \beta )\right) d\alpha d\beta \\
&& \left[ f''(Z(s)) D_{\alpha , \beta }Z(s) +
f'''(Z(s))D_{\alpha }Z(s) D_{\beta }Z(s)\right]  \\
&:=& C^{1} _{\varepsilon } + C^{2} _{\varepsilon }.
\end{eqnarray*}
We can write

\begin{eqnarray*}
C_{\varepsilon }&=& \frac{2}{\varepsilon }d(H) ^{2} \int _{0} ^{t}
f''(Z(s))ds \int _{0} ^{s} du \int _{s} ^{s+\varepsilon } dv \left(
\int _{0} ^{u\wedge v} \frac{\partial K^{H'} } {\partial u}(u, \alpha
)\frac{\partial K^{H'} } {\partial v}(v, \alpha )d\alpha \right) ^{2}\\
&=& \frac{2}{\varepsilon }d(H)^{2} (H'(2H'-1))^{2}  \int _{0} ^{t}
f''(Z(s))ds \int
_{0} ^{s}du\int _{s} ^{s+\varepsilon } dv \vert u-v\vert ^{2H-2} \\
&=& 2d(H)^{2} (H'(2H'-1))^{2} \int_{0}^{t} dv\int _{0} ^{v} du \vert
u-v\vert ^{2H-2}\left(  \frac{1}{\varepsilon } \int _{(v-\varepsilon
)\vee u} ^{v} f''(Z(s)) ds \right)
\end{eqnarray*}
and then clearly
\begin{equation}
\label{c1} C_{\varepsilon }^{1}  \to \frac{d(H)^{2} (H'(2H'-1))^{2}
} {2H-1} \int _{0}^{t} f''(Z(v) ) v^{2H-1} dv .
\end{equation}
in $L^{1}(\Omega )$ as $\varepsilon \to 0$.

The term denoted by $C_{\varepsilon }^{2} $ can be handled in the
following way:
\begin{eqnarray*}
C_{\varepsilon }^{2} &=& \frac{1}{\varepsilon } \int _{0}^{t}ds
f'''(Z(s)) \int _{0} ^{t} \int _{0} ^{t}D_{\alpha }Z(s) D_{\beta
}Z(s) \left( f_{s+\varepsilon }(\alpha , \beta ) -f_{s }(\alpha ,
\beta )\right) d\alpha d\beta \\
&=& \frac{4}{\varepsilon } d(H) ^{2} \int _{0}^{t}ds f'''(Z(s))
\int_{0} ^{s} \int _{0} ^{s} d\alpha d\beta \left( f_{s+\varepsilon
}(\alpha ,
\beta ) -f_{s }(\alpha , \beta )\right) \\
&& \left( \int _{0}^{t} dB(y_{1}) \int _{\alpha \vee y_{1}}^{s}
\frac{\partial K^{H'} }{\partial u} (u, \alpha )\frac{\partial K^{H'} }{\partial
u} (u, y_{1} )du \right)  \left( \int _{0}^{t} dB(y_{2}) \int
_{\beta \vee y_{2}} ^{s} \frac{\partial K^{H'} }{\partial v} (v, \beta
)\frac{\partial K^{H'} }{\partial v} (v, y_{2} )dv \right)
\end{eqnarray*}

\vskip0.5cm

{\bf The  case $f(x)=x^{3}$. }  In this case
\begin{eqnarray*}
C_{\varepsilon }^{2} &=& \frac{24}{\varepsilon } d(H)^{2} \int
_{0}^{t}ds
 \int_{0} ^{s}
\int _{0} ^{s} d\alpha d\beta \left( f_{s+\varepsilon }(\alpha ,
\beta ) -f_{s }(\alpha , \beta )\right) \\
&&\int _{0}^{t} \int _{0}^{t} dB(y_{1})dB(y_{2}) ) \int _{\alpha
\vee y_{1}}^{s} \frac{\partial K^{H'} }{\partial u} (u, \alpha
)\frac{\partial K^{H'} }{\partial u} (u, y_{1} )du
 \int
_{\beta \vee y_{2}} ^{s} \frac{\partial K^{H'} }{\partial v} (v, \beta
)\frac{\partial K^{H'} }{\partial v} (v, y_{2} )dv \\
&&+\frac{24}{\varepsilon } \int _{0}^{t}ds
 \int_{0} ^{s}
\int _{0} ^{s} d\alpha d\beta \left( f_{s+\varepsilon }(\alpha ,
\beta ) -f_{s }(\alpha , \beta )\right) \\
&&\int _{0} ^{s} dy_{1} \int _{\alpha \vee y_{1}}^{s} \frac{\partial
K}{\partial u} (u, \alpha )\frac{\partial K^{H'} }{\partial u} (u, y_{1}
)du
 \int
_{\beta \vee y_{2}} ^{s} \frac{\partial K^{H'} }{\partial v} (v, \beta
)\frac{\partial K^{H'} }{\partial v} (v, y_{2} )dv
\end{eqnarray*}
and by the same type of calculations as above we can prove that as $\varepsilon \to 0$
\begin{eqnarray}
C_{\varepsilon }^{2}&\to & 24 d(H)  ^{3} (H'(2H'-1)) ^{2} \int
_{0}^{t} \int _{0}^{t} dB(y_{1})dB(y_{2}) ) \int _{y_{1}}^{t}du'
\int _{y_{1}}^{u'} du\int
_{y_{2}} ^{u'} dv\\
&&\vert u-u'\vert ^{2H'-2} \vert v-u'\vert ^{2H'-2} \frac{\partial
K}{\partial u} (u, y_{1} )\frac{\partial K^{H'} }{\partial
v} (v, y_{2} ) \nonumber \\
&&+24d(H)^{3} (H'(2H'-1)) ^{3}\int _{0} ^{t} du' \int _{0} ^{u'}
\int _{0}^{u'} dudv\nonumber\\
&&\vert u-u'\vert ^{2H'-2} \vert v-u'\vert ^{2H'-2}
\vert u-v\vert ^{2H'-2}.
\end{eqnarray}

\vskip0.5cm

We can summarize
\begin{theorem}
We have
\begin{eqnarray*}
&& Z(t) ^{2} =2\int_{0}^{t} Z(s) \delta Z(s) + t^{2H} \nonumber \\
&+& \frac{4(2H-1) }{H+1} \int _{0} ^{t} \int
 _{0} ^{t} dB(y_{1})
 dB(y_{2})\int _{y_{2} }^{t}dv \int_{y_{1}} ^{u}du \frac{\partial K^{H'} } {\partial u}
(u,y_{1}) \frac{\partial K^{H'} } {\partial u} (u,y_{2}) \vert u-v\vert
^{2H'-2}.\label{itox2}
\end{eqnarray*}
and
\begin{eqnarray*}
&&Z(t)^{3} = 3 \int _{0}^{t} Z(s) ^{2} \delta Z(s) + 24\frac{ ( d(H) H'(2H'-1))^{2}}{2H-1} \int_{0}^{t} Z(s) s^{2H-1}ds\nonumber  \\
&& +24 d(H)^{2} H'(2H'-1)  \int
_{0}^{t} \int _{0} ^{t} dB(y_{1})
 dB(y_{2})\\
&&\int _{y_{2}} ^{t}dv\int _{y_{1}}^{v} duZ_{v}\frac{\partial K^{H'} }
{\partial u} (u, y_{1}) \frac{\partial K^{H'} } {\partial v} (v,y_{2})
\vert u-v\vert ^{ 2H'-2 }\nonumber \\
&& +72 d(H)^{3}( H'(2H'-1))^{2}  \int _{0}^{t}
\int _{0} ^{t} dB(y_{1})
 dB(y_{2})  \\
 && \int _{y_{2}} ^{t}dv \int _{y_{1}}^{v}du' \int _{y_{1}}^{v}
 du \vert u-u'\vert ^{2H'-2} \vert u-v\vert ^{2H'-2} \frac{\partial K^{H'} } {\partial u'} (u',
y_{1}) \frac{\partial K^{H'} } {\partial v} (v,y_{2}) \nonumber \\
&& +24d(H)^{3}( H'(2H'-1))^{3} \int_{0}^{t} dv \int _{0} ^{v} \int_{0}^{v} du'du
\vert u-v\vert ^{2H'-2}\vert u'-v\vert ^{2H'-2}\vert u-u'\vert
^{2H'-2}.
\end{eqnarray*}

\end{theorem}

\vskip0.5cm

\begin{remark}
One can note the appearance of a term involving $f'''$ in the
expression of the summand $C_{\varepsilon}$. Therefore one cannot
hope to have It\^o's formulas that end with a second derivative
term.
\end{remark}

 \end{document}